\documentclass[11pt,a4paper]{article}
\usepackage{epsfig}
\usepackage{graphics}
\usepackage{amsmath,amssymb}
\usepackage{subfigure}
\usepackage{graphicx}
\usepackage{epstopdf}
\usepackage{hyperref}
\usepackage{pdfsync}
\usepackage{color}

\numberwithin{equation}{section}
\newtheorem{thm}{Theorem}
\newtheorem{prop}{Proposition}
\newtheorem{cor}{Corollary}

\def \beqn {\begin{eqnarray}}
\def \beqnn {\begin{eqnarray*}}
\def \eeqn {\end{eqnarray}}
\def \eeqnn {\end{eqnarray*}}
\def \GG {\mathcal{G}}
\def \FF {\mathcal{F}}
\def \HH {\mathcal{H}}

\def \no {\Arrowvert}

\def \R {\mathbb{R}}

\def \P  {\mathbb{P}} 
\def \R {\mathbb{R}}
\def \E {\mathbb{E}}
\newtheorem{lemma}{Lemma}

\title{Noisy classification with boundary assumptions}
\author{S\'ebastien Loustau\footnote{Universit\'e d'Angers, LAREMA, loustau@math.univ-angers.fr} \ and Cl\'ement Marteau\footnote{Institut de Math\'ematiques de Toulouse, marteau@math.univ-toulouse.fr}}
\date{}

\begin{document}
\maketitle

\begin{abstract}
We address the problem of classification when data are collected from two samples with measurement errors. This problem turns to be an inverse problem and requires a specific treatment. In this context, we investigate the minimax rates of convergence using both a margin assumption, and a smoothness condition on the boundary of the set associated to the Bayes classifier. We establish lower and upper bounds (based on a deconvolution classifier) on these rates. 
\end{abstract}

\section{Introduction}

Assume that we have at our disposal two \textit{noisy} learning samples $\mathcal{S}_1=(Z_1^{(1)},\dots,Z_n^{(1)})$ and $\mathcal{S}_2=(Z_1^{(2)},\dots,Z_n^{(2)})$ satisfying:
\begin{equation}
Z_i^{(1)}=X_i^{(1)}+\epsilon^{(1)}_i, \ \forall i\in \lbrace 1,\dots, n\rbrace, \ \mathrm{and} \ Z_j^{(2)}=X_j^{(2)} + \epsilon^{(2)}_j, \ \forall j\in \lbrace 1,\dots, m \rbrace,
\label{eq:observation}
\end{equation}
where the $X_i^{(1)}$ (resp. $X^{(2)}_j$) denote independent identically distributed (i.i.d) random variables from unknown distribution $F_1$ (resp. $F_2$) and the $\epsilon_i^{(k)}$ denote random errors, independent of the $X_i^{(k)}$. For the sake of simplicity, $F_1$ (resp. $F_2$) admits a density $f$ (resp. $g$) w.r.t. the Lebesgue measure on $\mathbb{R}^d$. Moreover, we assume that the random errors are i.i.d. with \textit{known} density $\eta$ with respect to the Lebesgue measure. 

The goal is to classify a new incoming observation $X$, assumed to have a density $f$ or $g$ (and independent of $\mathcal{S}_1$ and $\mathcal{S}_2$). In other words, one wants to determine whether the density of $X$ is $f$ or $g$. Remark that in this model, two independent random sources are involved. The first one corresponds to the fluctuations of the variable of interest which is governed by the distribution $F_1$ or $F_2$ following the corresponding label. The second one corresponds to measurement errors (or imprecisions) during the data collection process. This problem corresponds to a nonparametric measurement error model, or errors-in-variables model (see the monograph of \cite{meister} for an introduction). We are in fact faced to the so-called discrimination with errors in variables problem.\\

The free noise case has already been widely studied in the literature. We refer for instance to  \cite{jaune} for a complete survey. When the variables $\epsilon_i^j$ are equal to zero in \eqref{eq:observation}, fast rates of convergence  were obtain for the first time in \cite{mammen}, using both a complexity and a margin assumption. Similar results where obtained in \cite{tsybakov2004,nedelec,AT} in slightly different settings. The previous papers were focused on empirical risk minimisation (ERM) algorithm and used margin assumptions (see Section \label{s:model} for more details). Similar conditions were investigated in the last decades for instance in \cite{localrademacher}, \cite{empimini} or \cite{kolt} among other.\\ 

Concerning the error-in-variables model of classification, few results have been published. Up to our knowledge, the only minimax result is \cite{les_pinkfloyd}, where minimax fast rates were obtained using a regularity assumption on the class of densities $f$ and $g$ in model \eqref{eq:observation} (see also \cite{inversestatlearning} in the general context of statistical learning). In a slightly different setting, we can however mention \cite{delaigle} where boundary estimation in a deconvolution framework was considered, or \cite{klemela} for a study of empirical risk minimization (ERM) algorithm in inverse problem models. Finally, \cite{koltgeometry} proposes an empirical minimization based on a deconvolution approach for the problem of estimating geometric characters of a multivariate distribution in the presence of noisy measurements.\\

The aim of this paper is to provide a classifier in the error-in-variables model and to study the related minimax performances in terms of fast rates. For this purpose, we will use two different kind of assumptions on the model: a margin assumption and a complexity assumption. The margin assumption will traduce the difficulty to discriminate an observation from another. It has been introduced for the first time in  \cite{mammen}. The second assumption concerns the regularity of the boundary of the set 
$$ G_K^\star = \left\lbrace x\in K \subset \R^d, f(x) \geq g(x) \right\rbrace,$$
where we restrict the study to a compact subset $K \subset \mathbb{R}^d$. It is widely known in classification that the decision set $G^*_K$ minimizes the so-called Bayes risk:
$$
R_K(G)=\frac{1}{2}\left(\int_{K\setminus G} f(x)dx+\int_G g(x)dx\right).
$$
Hence, the construction of a good classifier is more or less related to provide a good estimation (in a sense which will be precised later on) of $G_K^\star$. Remark that, contrary to \cite{AT} or \cite{les_pinkfloyd}, minimax results are investigated with no restriction over the regularity of $f$ and $g$.\\

The structure of this paper is as follows. In Section \ref{s:model}, we present in detail our model and assumptions. Then, we construct a deconvolution classifier. Lower bounds are provided in Section \ref{s:lowerbound}. The performances of our classifier are studied in Section \ref{s:upperbound}. Section \ref{s:conclusion} concludes the paper whereas Section \ref{s:preproof} proposes to highlight the main ideas used in the proofs. Section \ref{s:proof} is dedicated to the proofs of the main results, whereas Section \ref{s:appendix} adds some useful materials about noisy empirical processes.

\section{A deconvolution classifier}
\label{s:model}
\subsection{Model}
In this paper, a classifier is related to a subset $G\subset\mathbb{R}^d$ which traduces some hint on the places where there may be a greater probability to find an observation having distribution $F_1$. In the sequel, we restrict our investigations to a compact set $K \subset \mathbb{R}^d$. Using a slight abuse of notation, a classifier will be denoted by a measurable subset of the observations $\hat G=\hat G(\mathcal{S}_1,\mathcal{S}_2)$. In other words, the new incoming observation $X$ will be associated to the first (resp. second) label if it belongs to the set $\hat G$ (resp. $K\setminus \hat G$). \\

In order to measure the performances of a given classifier $\hat G \subset K$, we will use the Bayes risk defined as:
$$ R(\hat G)= \frac{1}{2} \left( \int_{K\setminus \hat G} f(x)dx + \int_{\hat G} g(x)dx \right).$$
The best possible classifier $G_K^\star$ then satisfies
$$ G_K^\star = \mathrm{arg} \min_{G \subset K} R_K(G) = \left\lbrace x\in K, f(x) \geq g(x) \right\rbrace,$$
where the infimum is taken over all possible subset of $K \subset \mathbb{R}^d$. In some sense, a good classifier should mimic (at least asymptotically) the behavior of $G_K^\star$. Hence, the excess risk
$$ R_K(\hat G) - R_K(G_K^*),$$
will be of first interest all along the paper.\\

For all $G\subset K$, using simple algebra, we get, 
\begin{eqnarray*}
R_K( G) - R_K(G_K^\star) & = & \frac{1}{2} \int_K (f-g)(x) (\mathbf{1}_{G_K^\star}(x) - \mathbf{1}_{G}(x))dx ,\\
& = & \frac{1}{2} \int_K |f(x)-g(x)| | \mathbf{1}_{G_K^\star}(x) - \mathbf{1}_{G}(x) |dx,
\end{eqnarray*}
since the product $(f(x)-g(x)).(\mathbf{1}_{G_K^\star}(x) - \mathbf{1}_{G}(x))$ is nonnegative for all $x\in K$. Then,
$$ R_K(\hat G) - R_K(G_K^\star) = \frac{1}{2} \int_K |f(x)-g(x)| \mathbf{1}_{G \Delta G_K^\star}(x)dx := \frac{1}{2}d_{f,g}(G,G_K^\star),$$
where for all $G_1,G_2 \subset K$, $G_1 \Delta G_2 = \lbrace K\setminus G_1 \cap G_2 \rbrace \cup \lbrace G_1 \cap K\setminus G_2 \rbrace$. The term $d_{f,g}$ is a pseudo distance on the subsets of $K$. The excess risk corresponds to a measure of the difference between $\hat G$ and the Bayes risk $G_K^{\star}$, where the symmetric difference is balanced by the value of $f-g$. This term is avoided when using for instance the pseudo-distance $d_\Delta$, defined as
$$ d_\Delta(G_1,G_2) = Q(G_1 \Delta G_2) \ \forall G_1,G_2 \subset K,$$
where $Q$ denotes the Lebesgue measure on $\R^d$.\\

\textsc{Remark 1.} The pseudo-distance $d_{f,g}$ is related to the densities $f$ and $g$ whereas $d_{\Delta}$ is entirely determined by the symmetric difference between $G_1$ and $G_2$. In some sense, $d_{f,g}$ is more related with the prediction task whereas $d_{\Delta}$ is a more related with a set estimation problem. In some favorable cases, i.e. when $Q(|f-g|\leq t_0)=0$ for some $t_0>0$, it is clear that these two pseudo-distance are equivalent (see the margin assumption below and Lemma 2 in \cite{mammen}). This particular case is known as the strong margin assumption case. \\

\textsc{Remark 2.} Our goal is to provide the best possible estimation of the set $G_K^\star$ from two noisy learning samples. From the prediction point of view, we are in fact interested in the estimation of the class of a new incoming observation X. We could also address the following problem: given a new noisy incoming observation, try to guess the corresponding label. These two problems are rather close but a precise comparison is beyond the scope of the present paper. We refer for instance to \cite{test2} where a similar problem was addressed in a goodness-of-fit testing framework, or to \cite{les_pinkfloyd} for a related discussion.\\

In this paper, our aim is to establish minimax rates of convergence for both $d_{f,g}$ and $d_\Delta$. In order to get these rates, we will need some assumptions on the model.

\subsection{Assumptions}

Following for instance \cite{mammen}, we will use two different conditions in order to obtain minimax rates of convergence. The first one is related to the behavior of the function $f-g$ at the boundary of $G_K^\star$. Recall that $Q$ denotes the Lebesgue measure on $\mathbb{R}^d$.\\
\\
\noindent
\textbf{Margin Assumption}: There exist constant $t_0,c_2\in\R_+$ and $\alpha\in\bar{\R}_+$ such that $\forall 0<t<t_0$,
\beqn
\label{ma}
Q\{x\in K:|f(x)-g(x)|\leq t\}\left\{
\begin{array}{l}
\leq c_2t^{\alpha}\mbox{ if }\alpha\in\R_+,\\
=0\mbox{ if }\alpha=+\infty.
\end{array}
\right.
\eeqn
This condition expresses the difficulty of distinguishing a distribution from another at the boundary of $G_K^\star$. It has been explicitly introduced for the first time in \cite{mammen}. The case $\alpha=+\infty$ corresponds to the best situation when $f-g$ does not hit or cross the frontier of the Bayes set. It is the so-called strong margin assumption. In this case, it is well-known from \cite{mammen} that $d_\Delta$ and $d_{f,g}$ are equivalent.  If $\alpha\in\R$, the most favorable cases corresponds to large values for $\alpha$: the distributions $F_1$ and $F_2$ are rather different of each side of $G_K^\star$. Small values for $\alpha$ correspond to more difficult situations where fast rates can not be expected.\\
\\
\textbf{Regularity Assumption.} The second condition is related to the complexity of the problem. Since we are dealing with a nonparametric set estimation problem, it seems natural to use an assumption on the regularity of the boundary of $G_K^\star$.  More precisely, we will deal with the family of boundary fragments on $K$. All along the paper, we state $K=[0,1]^d$ without loss of generality. A set $G\subset [0,1]^d$ belongs to a class of boundary fragments (see \cite{korostelevtsybakov}) if there exists $b:[0,1]^{d-1}\to [0,1]$ such that:
\beqnn
G=\{x=(x_1,\ldots x_d)\in [0,1]^d \ :x_d\leq b(x_1,\ldots,x_{d-1}) \}:=G_b.
\eeqnn
For given $\gamma,L>0$ the class of H\"older boundary fragments is then defined as:
\beqn
\label{hbfset}
\GG(\gamma,L)=\{G_b,b\in\Sigma(\gamma,L)\},
\eeqn
where $\Sigma(\gamma,L)$ is the class of isotropic H\"older continuous functions $b(x_1,\ldots,x_{d-1})$ having continuous partial derivatives up to order $\lfloor \gamma \rfloor$, the maximal integer strictly less than $\gamma$ and such that:
\beqnn
|b(y)-p_{b,x}(y)|\leq L|x-y|^\gamma,\forall x,y\in \mathbb{R}^{d-1},
\eeqnn
where $p_{b,x}$ is the Taylor polynomial of $b$ at order $\lfloor \gamma \rfloor$ at point $x$. \\
\\
In the sequel, we restrict the class $\mathcal{G}$ of possible candidate sets $G\subset K$ for which both the margin and H\"older boundary fragment assumptions are satisfied. It requires, in turn, restrictions on the class $\mathcal{F}$ of possible density couple $(f,g)$. Our result are given in a minimax framework over the following class $\mathcal{F}$. For positive constants $\gamma, L, c_2,t_0, \alpha$ and $c_1$, the class $\mathcal{F}$ is defined as:
\beqn
\label{classF}
\mathcal{F}(\alpha,\gamma)=\{(f,g): f \mbox{ and } g\mbox{ are densities w.r.t. the Lebesgue measure},\no f\no_\infty\vee\no g\no_\infty\leq c_1,\nonumber\\
\{x\in K:f(x)\geq g(x)\}\in\mathcal{G}(\gamma,L)\mbox{ and \eqref{ma} holds for $\alpha\in\bar{\R}$}
\}.
\eeqn

In the free-noise case, \cite{mammen} has proved that the rates in $d_\Delta$ and $d_{f,g}$ can be completely characterized by both margin and boundary fragment assumptions. Remark that alternative hypotheses can be set on the model. For instance, \cite{AT} or \cite{les_pinkfloyd} deal with a plug-in type assumption on the regularity of $f-g$.\\

The last hypothesis that we will introduce on the model concerns the measurement errors. Indeed, in the model (\ref{eq:observation}), the density of the $Z_i^{(1)}$ (rep. $Z_i^{(2)}$) is nor $f$ (resp. $g$) but rather $f*\eta$ (resp. $g*\eta$), where $*$ denotes the convolution product between two functions and $\eta$ the density of the $\epsilon_i^{(j)}$ w.r.t. the Lebesgue measure. Contrary to the free-noise case, the $X_i^{(j)}$ are indirectly observed: we are faced to an inverse (deconvolution) problem. 

Inverse problems have been widely investigated in the statistical literature. We mention for instance \cite{meister} or \cite{lechef} for a general review of existing models and related results. In an estimation or testing framework, inverse problems are known for providing slower rates than in the direct cases. This can be explained by the loss of information related to the regularization of the operator. The behavior of the noise density $\eta$ is hence of first importance if one want to evaluate this decay. In particular, we will see that we can take advantage of the shape of the Fourier transform of the noise in order to provide a precise description of the minimax rates in this setting. This is the purpose of the following assumption. \\

\noindent
\textbf{Noise Assumption}: \textit{There exists $\beta=(\beta_1,\dots,\beta_d)'\in \R_+^d$ such that for all $i\in \lbrace 1,\dots, d \rbrace$, $\beta_i>1/2$,
$$ \left| \mathcal{F}[\eta_i](t) \right| \sim |t|^{-\beta_i}, \ \mathrm{and} \ \left| \mathcal{F}'[\eta_i](t) \right| \leq C|t|^{-\beta_i} \ \mathrm{as} \ t\to +\infty,$$
where $\eta=\Pi_{i=1}^d\eta_i$ and $\mathcal{F}[\eta_i]$ denotes the Fourier transform of $\eta_i$. Moreover, we assume that $\mathcal{F}[\eta_i](t) \not = 0$ for all $t\in \R$ and $i\in \lbrace 1,\dots, d \rbrace$.}\\

Note that the hypothesis $\eta=\Pi_{i=1}^d \eta_i$ corresponds to non-degenerated random errors $\epsilon$, whose coordinates are independent. This assumption could be relaxed as in $\cite{comtelacour}$ since we only need an assumption over the asymptotic behavior of $\mathcal{F}[\eta]$. In the literature, the noise assumption \textbf{(NA)} corresponds to \textit{ordinary smooth} or \textit{mildly ill-posed} inverse problem. The parameters $\beta_i$ describe the difficulty of the related problem. Higher is $\beta_i$, smoother are $f*\eta$ and $g*\eta$ in the direction $i$. As a result, harder becomes the classification problem. In the sequel, we show that these coefficients play a crucial role in the expression of the minimax rates of convergence. For the sake of concision, we will not consider \textit{severely ill-posed} problems, i.e. corresponding to exponentially decreasing Fourier transform. However, simple applications of the main result of this paper lead to the study of this particular case.\\

\subsection{The classifier}

We are now ready to propose a classifier in such a context. Our method is based on the empirical risk minimization (ERM) method. The main idea is to construct an estimator for the Bayes risk associated to each candidate $G$, and then to select the one associated to the lowest value. In the free-noise case, i.e. when data are observed without measurement errors, \cite{mammen} have used the risk estimator $R_{n,m}(.)$ defined as
$$ R_{n,m}(G) = \frac{1}{2} \left[ \frac{1}{n} \sum_{i=1}^n \mathbf{1}_{\lbrace X_i^{(1)} \in K\setminus G \rbrace} + \frac{1}{m} \sum_{i=1}^m \mathbf{1}_{\lbrace X_i^{(2)} \in G \rbrace} \right], \ \forall G\subset K.$$
In particular, it is easy to see that for all $G\subset K$, $ R_{n,m}(G)$ is an unbiased and consistent estimator of $R_K(G)$. When dealing with an error-in-variables model, the methodology is completely different. Indeed, for all $i\in \lbrace 1,\dots, n \rbrace$, we get for instance
$$ \mathbb{E}\left[ \mathbf{1}_{\lbrace Z_i^{(1)} \in K\setminus G \rbrace} \right] = \int_{K\setminus G} f*\eta(x)dx.$$
In such a situation, $ R_{n,m}(G)$ is nor an unbiased neither a consistent estimator of the risk $R_K(G)$. We are faced to an inverse (deconvolution) problem. In order to get round of this problem, we will propose a deconvolution ERM algorithm. This algorithm is heavily related to the properties of deconvolution kernel (see for instance \cite{fan} or \cite{meister}).  \\


Let $\mathcal{K}=\prod_{j=1}^d \mathcal{K}_j:\R^d \to \R$ be a $d$-dimensional kernel defined as the product of $d$ unidimensional kernels $\mathcal{K}_j$ (i.e. functions $\mathcal{K}_j : \mathbb{R} \rightarrow \mathbb{R}$ satisfying $\int \mathcal{K}_j =1$). The properties of $\mathcal{K}$ leading to satisfying upper bounds will be made precise later on. Then, if we denote by $\lambda=(\lambda_1,\dots,\lambda_d)$ a set of (positive) bandwidths and by $\FF[\cdot]$ the Fourier transform, we define the deconvolution kernel $\mathcal{K}_\eta$ as
\begin{eqnarray}
\mathcal{K}_{\eta} & : & \R^d \to \R \nonumber \\
& & t \mapsto \mathcal{K}_\eta(t) = \FF^{-1}\left[ \frac{\FF[\mathcal{K}](\cdot)}{\FF[\eta](\cdot/\lambda)}\right](t),
\label{dk}
\end{eqnarray}
provided that $\mathcal{K}$ (resp. $\eta$) belongs to $L_2(\R^d)$ and admits a Fourier transform. Note that in the sequel, for the sake of concision, we note, for any $\lambda\in\R^d_+$, $z,x\in\R^d$:
$$
\frac{1}{\lambda} \mathcal{K}_\eta \left(\frac{z-x}{\lambda}\right):=\frac{1}{\lambda_1\cdots\lambda_d}\mathcal{K}_\eta \left(\frac{z_1-x_1}{\lambda_1},\ldots,\frac{z_d-x_d}{\lambda_d}\right).
$$
In this context, for all $G\subset K$, the risk $R_K(G)$ can be estimated by
\beqnn
R^{\lambda}_{n,m}(G)=\frac{1}{2} \left[\frac{1}{n}\sum_{i=1}^nh_{K\setminus  G,\lambda}(Z_i^{(1)})+\frac{1}{m}\sum_{i=1}^mh_{G,\lambda}(Z_i^{(2)})\right],
\label{eq:empriskbruit1}
\eeqnn
where for a given $z\in \R^d$:
\beqn
\label{hG}
h_{G,\lambda}(z)= \int_{G}\frac{1}{\lambda} \mathcal{K}_\eta \left(\frac{z-x}{\lambda} \right)dx.
\eeqn
For all $G\subset K$, the function $h_{G,\lambda}$ more or less plays the role of an indicator function. In particular, for all $i\in \lbrace 1,\dots, n \rbrace$, we get for instance 
\begin{equation}
\mathbb{E}[h_{G,\lambda}(Z_i^{(1)}) / X_i^{(1)}] = \mathcal{K}_\lambda*\mathbf{1}_{\lbrace .\in G\rbrace}(X_i^{(1)}), \ \forall G\subset K,
\label{eq:biais1}
\end{equation}
where $\mathcal{K}_\lambda*\mathbf{1}_{\lbrace .\in G\rbrace}(x)$ denotes the convolution between $\mathcal{K}$ and the indicator function $\mathbf{1}_{\lbrace .\in G\rbrace}$ at a point $x\in \mathbb{R}$. This term can then be viewed as a smoothed indicator on the set $G$. Remark that due to (\ref{eq:biais1}), the estimator (\ref{eq:empriskbruit1}) will be biased. Indeed, for all $G\subset K$
$$ \mathbb{E} R_{n,m}(G) = \int_{\mathbb{R}^d} f(x)  \mathcal{K}_\lambda*\mathbf{1}_{\lbrace .\in K/G\rbrace}(x) + \int_{\mathbb{R}^d} g(x)  \mathcal{K}_\lambda*\mathbf{1}_{\lbrace .\in G\rbrace}(x):=R^\lambda(G) \not = R(G).$$ 
The control of the related bias will be one of the main difficulty to establish minimax rates of convergence for this estimator.\\

In the following, we study ERM estimators defined as:
\beqn
\label{derm2}
\hat{G}_{n,m}^\lambda=\arg\min_{G\in\mathcal{G}(\gamma,L)} R_{n,m}^\lambda(G),
\eeqn
where $\mathcal{G}(\alpha,\gamma)$ is defined in \eqref{hbfset} and $\lambda=(\lambda_1,\ldots ,\lambda_d)\in\R_+^{d}$ is a parameter that has to be chosen explicitly.



\section{Lower bound}
\label{s:lowerbound}

Theorem 1 states lower bounds for the minimax risks over the class $\mathcal{F}(\alpha,\gamma)$ defined in \eqref{classF}. The proof is postponed to Section \ref{section:proofs}.

\begin{thm}
\label{thm:lbboundary}
Let $K=[0,1]^d$ and $\mathcal{F}(\alpha,\gamma)$ defined in \eqref{classF}. Suppose that the noise assumption is satisfied for some $\beta$. Then we have:
$$\liminf_{n \to +\infty} \inf_{\hat G_{n,m}} \sup_{(f,g)\in \mathcal{F}(\alpha,\gamma)} (n\wedge m)^{\tau_d(\alpha,\beta,\gamma)} \mathbb{E} d_{\square}(\hat G_{n,m},G_K^{\star}) > 0,$$
where the infinimum is taken over all possible estimators of the set $G_K^{\star}$ and
$$ \tau_d(\alpha,\beta,\gamma) =\left\{
\begin{array}{l}
 \frac{\displaystyle \gamma\alpha }{\displaystyle \gamma(2+\alpha)+(d-1)\alpha+ 2\alpha\sum_{i=1}^{d-1}\beta_i +2\alpha\beta_d \gamma} \mbox{ for }d_{\square}=d_\Delta\\
 \\
 \frac{\displaystyle \gamma(\alpha+1)}{\displaystyle \gamma(2+\alpha)+(d-1)\alpha+ 2\alpha\sum_{i=1}^{d-1}\beta_i +2\alpha\beta_d \gamma} \mbox{ for } d_{\square}=d_{f,g}.
\end{array}\right.
$$
\end{thm}

\textsc{Remark 4.} We obtain exactly the same lower bounds as \cite{mammen} in the direct case, which yet corresponds to the situation where $\beta_j=0$ for all $j\in \lbrace 1, \dots , d \rbrace$. In this particular framework, the minimax rate of convergence mainly depends on $\gamma$ and $\alpha$. The coefficient $\gamma$ corresponds to the regularity of the boundary of $G_K^\star$. Greater is $\gamma$, easier is the estimation. The term $\alpha$ is related to the margin assumption.\\

\textsc{Remark 5.} In the presence of noise in the variables, the rates obtained in Theorem \ref{thm:lbboundary} are slower. The price to pay is an additional term of the form 
$$ 2\alpha\left[\sum_{i=1}^{d-1}\beta_i +\beta_d \gamma\right].$$
This term clearly connects the difficulty of the problem to the values of the coefficients $\beta_1,\dots,\beta_d$. Moreover, the above expression highlights a connection between the margin parameter and the ill-posedness. The role of the margin parameter over the inverse problem can be summarized as follows. Higher is the margin, higher is the price to pay for a given degree of ill-posedness. When the margin parameter is small, the problem is difficult at the boundary of $G^\star_K$ and we can only expect a non-sharp estimation of $G^\star_K$. In this case it is not significantly worst to add noise. On the contrary, for large margin parameter, there is nice hope to give a sharp estimation of $G^\star_K$ and then perturb the inputs variables have strong consequences in the performances. \\


\textsc{Remark 6.} In the above expression, the first $d-1$ components of $\epsilon$ do not have the same impact as the last (vertical) component. This is due to the fact that we consider boundary fragments with a given regularity $\gamma$. This regularity is expressed in a H\"older space of functions defined on the $d-1$ first directions.\\

\textsc{Remark 7.} Finally, we can compare the lower bound of Theorem \ref{thm:lbboundary} with the previous lower bound stated in \cite{les_pinkfloyd} under plug-in type conditions. The main novelty here is that no restriction on $\alpha\in\bar{\R}$ is necessary to get the lower bound. It could be explain as follows. Coarsely, the case $\alpha=+\infty$ cannot be treated in \cite{les_pinkfloyd} since the minimax approach is performed over a class of densities with H\"older regularity. If the strong margin assumption holds,  $f-g$ is not continuous at the boundary. Moreover, Theorem 1 holds for arbitrary values for $\alpha$. Since we do not suppose any assumption for the regularity of the densities $f$ and $g$, the construction of the lower bound is easier. In particular, we can take advantage of the noise assumption and mix standard arguments from lower bounds in classification (see \cite{audibert2004} and \cite{mammen}) and inverse problems (see \cite{Butucea}). 


\section{Upper bounds}
\label{s:upperbound}

\subsection{A preliminary result}
For the sake of concision, in this section we propose to restrict the set $\mathcal{G}$ to $\mathcal{G}(\gamma,L)$, where all possible regularities $\gamma$ satisfy $\gamma>d-1$. It allows us to control the bracketing entropy of $\mathcal{G}(\gamma,L)$ with a parameter $\rho=\frac{d-1}{\gamma}<1$. We may also consider more general classes of candidates $\mathcal{G}$, with given entropy rates. This extension is presented in Section \ref{s:preproof} using empirical processes theory in a more general framework.

In this section, we are interested in the performances of the estimator:
\beqn
\hat G_{n,m} = \mathrm{arg} \min_{G \in \mathcal{G}(\gamma,L)} R_{n,m}^\lambda(G),
\label{eq:sieve2}
\eeqn
where $\mathcal{G}(\gamma,L)$ is defined in \eqref{hbfset} for $\gamma>d-1$. Nevertheless, one may also define our ERM estimator for $\gamma\leq d-1$ by considering a network in a practical purpose, without significant change in the following results. We will also assume for clarity throughout this section that $n=m$. In order to get round of the assumption on both the shape of the noise and the boundary of $G_K^\star$, we will introduce some constraints on the kernel $\mathcal{K}$.\\

\noindent
\textbf{Kernel Assumption}
\begin{itemize}
\item[\textbf{(K1)}] The kernel $\mathcal{K}$ is such that the associated deconvolution kernel satisfies
$$ \sup_{t\in \mathbb{R}^d} \left| \mathcal{F}[\mathcal{K}_\eta] (t) \right| \leq C \prod_{i=1}^d \lambda_i^{-\beta_i}, \mathrm{and} \ \| \mathcal{K}_\eta \|^2 \leq C\prod_{i=1}^d \lambda_i^{-2\beta_i}. $$
\end{itemize}

The assumption \textbf{(K1)} is necessary to control the variance of our classifier. It is satisfied for instance if the Fourier transform of $\mathcal{K}$ is bounded and compactly supported.\\

The following theorem provides a control of the expectation of the excess risk by a bias-variance decomposition, up to a residual term depending on the choice of the bandwidth $\lambda$.
 
\begin{thm}
\label{prop:aminima}
Let $\hat G_{n}$ the set introduced in \eqref{eq:sieve2} where $\gamma>d-1$ and $n=m$. Suppose that the noise assumption is satisfied and consider a kernel $\mathcal{K}_\eta$ defined as in \eqref{dk} satisfying \textbf{(K1)}. Then we have
$$
\E d_{f,g}(\hat G_{n,m},G^*)\leq C\inf_{\lambda\in\R^d_+}\left[\left(\frac{\Pi_{i=1}^d\lambda_i^{-\beta_i}}{\sqrt{n}}\right)^{\frac{2\gamma(\alpha+1)}{\gamma(\alpha+2)+(d-1)\alpha}} +\sup_{G\in\mathcal{G}}\left(R_K-R_K^\lambda\right)( G)+\sum_{i=1}^d(n\lambda_i)^{-\frac{\gamma}{\gamma+d-1}}\right],
$$
where $C$ is a positive constant and $(R_K-R_K^\lambda )( G)=R(G)-R^\lambda(G)$ for all $G\in \mathcal{G}$.
\end{thm} 
This result highlights a bias-variance decomposition of the excess risk. The proof is presented in Section \ref{s:proof}. The main ingredient of the proof is a study of the increments of a noisy empirical process, indexed by a set of functions which depends on the regularization parameter $\lambda>0$. At this step, some remarks are necessary:\\

\textsc{Remark 9.} The variance term is obtained thanks to extensions of the empirical process machinery and the peeling technique introduced by \cite{vdg} in the direct case (see Section \ref{s:appendix} for details). This term is related to the regularity of the distribution function $\eta$ in the noise assumption. We can see coarsely that the price to pay for the inverse problem in the variance is summarized in the term $\Pi_{i=1}^d\lambda_i^{-\beta_i}$. Note that in the direct case, \cite{mammen} has already stated fast rates of the form $n^{-\frac{\gamma(\alpha+1)}{\gamma(\alpha+2)+(d-1)\alpha}}$, which corresponds to $\beta=0$ in Theorem \ref{prop:aminima}.\\

\textsc{Remark 10.} The second term in Theorem \ref{prop:aminima} is a bias term due to the estimation of the true risk by a biased empirical risk. When dealing with a deconvolution ERM, the algebra is rather different. We have to provide a precise control of the bias of the ERM, namely the quantity
$$
R(G)-\E R_n^\lambda(G)=\int (f-g)\left(\mathbf{1}_{G^C}-\mathcal{K}_\lambda*\mathbf{1}_{G^C}\right)dQ.
$$
This term has to be controlled carefully to get minimax results. This is the focus of the next paragraph.\\

\textsc{Remark 11.} Finally, the last term in the upper bound is a residual term since we can see coarsely that
$$
 (n\lambda_i)^{-\frac{\gamma}{\gamma+d-1}}\leq \left(\frac{\Pi_{i=1}^d\lambda_i^{-\beta_i}}{\sqrt{n}}\right)^{\frac{2\gamma(\alpha+1)}{\gamma(\alpha+2)+(d-1)\alpha}},
 $$
 provided that $\lambda_i\to 0$ not too fast, for all $i\in\{1,\ldots, d\}$. \\


\subsection{Control of the bias and related rates of convergence}
The aim of this part is to investigate different available ways in order to control the bias. It is important to note that in a previous paper dedicated to plug-in type conditions, \cite{les_pinkfloyd} provides a simple way to control the bias term. Indeed, under a H\"older regularity condition over the function $f-g$, we can bound the bias term as follows:
\begin{lemma}[Loustau and Marteau \cite{les_pinkfloyd}]
\label{pflemma}
Suppose $f-g\in\Sigma(\gamma,L)$, the isotropic H\"older class of functions over $\R^d$. Suppose that the kernel $\mathcal{K}$ is of order $\gamma$. Then, we have
$$
\sup_{G\subset K}(R^\lambda_K-R_K)(G)\leq C\sum_{i=1}^d\lambda_i^{\gamma}.
$$
\end{lemma}
The proof is straightforward since in this case, we can write:
\beqnn
R(G)-\E R_n^\lambda(G)=\int \mathbf{1}_{G^C}[(f-g)-\mathcal{K}_\lambda*(f-g)]dQ.
\eeqnn
Then, the control of the bias term is reduced to the control of the bias term in standard nonparametric density estimation, which gives (see for instance \cite{booktsybakov}):
 $$\sup_{x_0\in\R^d}\left|(f-g)(x_0)-\mathcal{K}_\lambda*(f-g)(x_0)\right|\leq \sum_{i=1}^d\lambda_i^{\gamma}.$$

Here, the problem is rather different since the regularity assumption deals with the boundary of $G^*_K$. It is well-known (see for instance \cite{korostelevtsybakov}) that a regularity with respect to the boundary of  a decision rule does not match with plug-in type conditions of Lemma \ref{pflemma}. The following result proposes an upper bound under boundary assumptions.

\begin{cor}
\label{cor:rough}
Let $\hat G_{n}$ the set introduced in \eqref{eq:sieve2} where $\gamma>d-1$. Suppose the noise assumption is satisfied and consider a kernel $\mathcal{K}_\eta$ defined as in \eqref{dk} satisfying \textbf{(K1)}. Suppose moreover that for any $j\in\{1,\ldots, d\}$, $\int_{\R^{d}}\left|\mathcal{K}(z)||z_j\right|dz<\infty$ and $\Pi_{j=1}^{d-1}\mathcal{K}_j$ has compact support. Then, there exists a positive constant $C$ such that
$$\E d_{\square}(\hat G_{n,m},G_K^*)\leq C n^{-\kappa_d(\alpha\beta,\gamma)},$$
where
$$ \kappa_d(\alpha,\beta,\gamma) =\left\{
\begin{array}{l}
 \frac{\displaystyle\gamma\alpha}{\displaystyle\gamma(\alpha+2) + (d-1)\alpha + 2\gamma(\alpha+1) \sum_{i=1}^d \beta_i}  \mbox{ for }d_{\square}=d_\Delta\\
 \\
 \frac{\displaystyle\gamma(\alpha+1)}{\displaystyle\gamma(\alpha+2) + (d-1)\alpha + 2\gamma(\alpha+1) \sum_{i=1}^d \beta_i} \mbox{ for } d_{\square}=d_{f,g}.
\end{array}\right.
$$
\end{cor}
Following Corollary \ref{cor:rough}, lower and upper bounds do not match. The prize to pay for the errors-in-variables model is summarized in the term $2\gamma(\alpha+1)\sum_{i=1}^d\beta_i$ whereas the lower bound proposes a smaller term $2\alpha\sum_{i=1}^{d-1}\beta_i+2\gamma\alpha\beta_d$. By the way, the corresponding error becomes negligible when $\gamma$ is close to $1$ and $\alpha\to\infty$. The proof is provided in Section \ref{s:proof} and uses Theorem \ref{prop:aminima} gathering with the following crude bound for the bias term:
$$
\sup_{G\in\mathcal{G}}(R^\lambda_K-R_K)(G)\leq C\sum_{i=1}^d\lambda_i.
$$
It is based on the following scheme. For all $G\subset K$, using Fubini, we have
\begin{eqnarray*}
\lefteqn{\int_{\mathbb{R}^d} (f-g)(x) \left( \mathcal{K}_\lambda*\mathbf{1}_G(x) - \mathbf{1}_G(x) \right) dx}\\
& = & \int_{\mathbb{R}^d} (f-g)(x) \left( \int_{z\in \mathbb{R}^2} \mathcal{K}(z) \left[ \mathbf{1}_G(x+\lambda z) - \mathbf{1}_G(x) \right]dz \right) dx\\
& = & \int_{\mathbb{R}^d} \mathcal{K}(z) \left( \int_{\mathbb{R}^d} (f-g)(x)  \left[ \mathbf{1}_G(x+\lambda z) - \mathbf{1}_G(x) \right]dx \right) dz.
\end{eqnarray*}
Since we do not have any conditions on the smoothness of $f-g$, the control of the bias reduces to the calculation of the Lebesgue measure between the sets $G$ and $G+\lambda z$, which appears to be of order $\sum_i \lambda_i$. Hence, we can not take advantage on the smoothness of the boundary. In Section \ref{s:conclusion} below, we discuss several tracks to attack this problem. It appeals to different tools (such as convexification, or additional regularity assumptions) which do not fit with the machinery of the present paper.\\

\section{Conclusion}
\label{s:conclusion}
Let us discuss the obtained results and highlight some open problems:
\begin{description}
\item[Comparison with \cite{mammen}] This paper can be seen as a generalization of the results of \cite{mammen} to the error-in-variables case. We highlight, in the presence of noise, fast rates of convergence which depends on the Fourier transform of the noise distribution $\eta$. The price to pay depends on the triplet $(\gamma,\alpha,\beta)$ related with the regularity, margin and noise assumptions.  

\item[Choice of $\lambda$ and model selection] The main drawback of the deconvolution ERM of this paper is the calibration of the bandwidths $\lambda$. Under isotropic assumptions over the shape of $f$ and $g$, we have shown that it is sufficient to choose only one bandwidth. An extension to the anisotropic case could be done easily (see for instance \cite{aoinq} in an unsupervised context). However, these calibrations are non-adaptive and depend on the smoothness assumptions. The data-driven choice of the bandwidth is a natural open problem. The bias variance trade-off to choose $\lambda$ is not the usual one in non-parametric statistics and a careful study of this problem is necessary. In this direction, we can mention the recent work of \cite{loustauchichi}. Moreover, this problem of adaptation is compounded with the model selection of $\mathcal{G}$.
\item[Adaptation to the operator] The deconvolution classifier proposed in this paper depends on the known density of the noise $\epsilon$. As a result, another issue would be to try to adapt to unknown error densities $\eta$. In this direction, it can be interesting to apply the same strategy, using for instance the estimator proposed in \cite{delaiglehallmeister} for deconvolution with repeated measurements. In the presence of repeated measurements, the model \eqref{eq:observation} becomes, for $i=1,\ldots, n$ and $ j=1,\ldots ,m$:
$$Z_{i,k}^{(1)}=X_i^{(1)}+\epsilon_{ik}, \,k=1,\ldots, N_i, \mbox{ and }  Z_{j,\ell}^{(2)}=X_j^{(2)}+\epsilon_{j\ell}, \,\ell=1,\ldots, M_j.
$$
In this case, the empirical risk associated to this problem can be written:
\beqnn
R^{\lambda}_{n,m}(G)=\frac{1}{2} \left[\frac{1}{nN}\sum_{i=1}^nh_{K\setminus G,\lambda}(Z_{i,1}^{(1)},\ldots,Z_{i,N_i}^{(1)})+\frac{1}{mM}\sum_{j=1}^mh_{G,\lambda}(Z_{j,1}^{(2)},\ldots,Z_{j,M_j}^{(2)})\right],
\label{eq:empriskbruit}
\eeqnn
where $N=\sum_{i=1}^nN_i$, $M=\sum_{j=1}^m M_j$, and for some $j\in\{1,\ldots,m\}$:
\beqn
\label{hGrm}
h_{G,\lambda}(Z_{j,2}^{(1)},\ldots,Z_{j,M_j}^{(2)})= \int_{G}w_j\sum_{\ell=1}^{M_j}\frac{1}{\lambda} \hat{L} \left(\frac{x-Z_{j,\ell}^{(2)}}{\lambda} \right)dQ(x),
\eeqn
where $w_j$ are weights satisfying $\sum_{j=1}^m w_jM_j=M$ and $\hat{L}$ is a ridged deconvolution estimator defined in \cite{delaiglehallmeister}. An interesting open problem is to use the same methodology presented in this paper to study the performances of the minimizer of the deconvolution ERM using \eqref{hGrm}. 

\item[Direct VS inverse problem] In this paper, we propose to study the rates of convergence to the Bayes of our deconvolution classifier in term of pseudo-distance $d_{f,g}$ and $d_{\Delta}$. However, as discussed in \textsc{Remark 2}, a direct approach seems to be tractable in some particular cases. A systematic study of the rates of convergence of standard ERM using noisy measurements could be interesting. For this purpose, we have to control the difference between the Bayes risk with respect to $X$ and the following  Bayes risk:
$$
R_{K}^{\eta}(G)=\frac{1}{2}\left[\int_{K\setminus G}f*\eta dQ+\int_G g*\eta dQ\right].
$$
Another interesting comparison could be done in the problem of predicting a new incoming noisy observation. This paper shows that deconvolution ERM are optimal to predict a new $X$ observation. Finding the more efficient rule to predict a new $Z$ observation is also of practical interest. To this end, a relationship between the margin assumption and the noise assumption has to be done, which appears to be a challenging open problem. 

\item[Minimax optimality remains an open problem] Lower and upper bounds of Theorem \ref{thm:lbboundary} and Corollary \ref{cor:rough} do not match and the question of minimax rates in such a setting remains an open problem. However, our intuition is the following: the lower bound is valid and on the contrary, the estimation method of this paper suffers from a lack of optimality.

Firstly, an investigation of the upper bound above indicates that an \textit{optimal} control of the bias would require an upper bound of order
$$ \left[ \left( \sum_{i=1}^{d-1} \lambda_i \right)^\gamma + \lambda_d \right]^{1+\frac{1}{\alpha}} \ \ \mathrm{instead} \  \mathrm{of} \ \ \sum_{i=1}^d \lambda_i.$$
Using simple algebra (in dimension 2 for the sake of convenience), we can re-write the bias as
\begin{eqnarray*}
\lefteqn{\int_{\mathbb{R}^d} (f-g)(x) \left( \mathcal{K}_\lambda*\mathbf{1}_G(x) - \mathbf{1}_G(x) \right) dx}\\
& = & \hspace{-0.4cm} \int_{\mathbb{R}^2} \mathcal{K}(z) \int_\mathbb{R} \left[\int_0^{(b(x_1+\lambda_1 z_1) - \lambda_2 z_2)_+} (f-g)(x) \mathbf{1}_{\lbrace x_1+\lambda_1 z_1 \in [0,1] \rbrace}dx_2 - \int_0^{b(x_1)} (f-g)(x) dx_2 \right] dx_1dz.
\end{eqnarray*}
The exponent $\gamma$ is related to the smoothness of the boundary of $b$. This smoothness properties could certainly be taken into account following some additional assumptions on the smoothness of $f-g$. Indeed, one may manage a Taylor expansion of $b$ and then $f-g$ in the $2^{nd}$ direction, up to some additional technical constraints. Concerning the exponent $1+1/\alpha$, one might take advantage of the behavior of $f-g$ in the neighborhood of $G_K^\star$, but this may require an extended version of the margin assumption. 

Another possibility is to use convex surrogate to deal with the bias term. Indeed, the difficulty to control the bias term seems to be related with the $0-1$ loss approach of this paper. Since we use the $0-1$ loss, the bias is upper bound by the Lebesgue measure between the sets $G$ and $G+\lambda z$, which is of order $\sum\lambda_i$ (see Corollary \ref{cor:rough} and the associated discussion). There is nice hope that a control of the bias term can be managed thanks to a smooth loss function, without any additional smoothness assumption. However, in this case, a precise study of the lower bound has to be performed and is out of the scope of the present paper. 

Finally, the lack of optimality of the deconvolution ERM of this paper could be explain as follows. In the estimation procedure, the idea is to estimate the true risk by using a deconvolution kernel estimator of the densities $f$ and $g$. As a result, the method seems to be strongly linked with plug-in type regularities for the Bayes decision rule $G^*_K$, which allows us to control easily the bias term in \cite{les_pinkfloyd}. Here, the regularity assumption is rather different and deals with the boundaries of $G^*_K$. That's why an optimal control of the bias term seems problematic. Another way to obtain the fast rates of the lower bound should be to use another estimator of the true risk to deal with errors-in-variables. In classification with smooth boundaries, the estimation task could be summarize as an estimation of the boundary. As a result, in the presence of noisy observations, we have to study the effect of the inverse problem on the boundary of the Bayes $G^*_K$ (that is with respect to the direct data). It could be a way to plug another estimator in the true risk, which allows to estimate optimally the boundary. 
\end{description}

\section{Preliminaries to the proofs}
\label{s:preproof}
In this section, we provide a more general point of view about the problem of classification with noisy inputs. To be more precise, we state a generalization of Theorem \ref{prop:aminima} to study any possible candidate set $\mathcal{G}$ in the ERM minimization, given its entropy rates.\\ 

More formaly, we suggest to state upper bounds for deconvolution ERM of the form:
\beqn
\label{generalsetestimation}
\hat{G}_n=\arg\min_{G\in\mathcal{G}}R_n^\lambda(G),
\eeqn
where $R_n^\lambda(\cdot)$ is the empirical risk defined in Section 1 and $\mathcal{G}$ is a set of possible candidates for $G^\star_K$. We want to study the rate of convergence of $\hat{G}_n$ to $G^\star_K$. This rate depends on the complexity of the class $\mathcal{G}$ in terms of $\delta-$entropy with bracketing. For $\delta>0$, the bracketing entropy of $\mathcal{G}$ with respect to some distance $d$ is denoted by $\mathcal{H}(\mathcal{G},d,\delta)$ and corresponds to the minimal number such that $\mathcal{N}_B(\delta)=\exp\left(\mathcal{H}(\mathcal{G},d,\delta)\right)$ is an integer and such that there exists pairs $(G_j,H_j)$, $j=1,\ldots, N_B(\delta)$ of subsets of $\mathcal{G}$ satisfying
\begin{itemize}
\item[(1)] $G_j\subset H_j$ for all $j\in \lbrace 1,\dots,  N_B(\delta)\rbrace$,
\item[(2)] $d(G_j,H_j)\leq \delta$ for all $j\in \lbrace 1,\dots,  N_B(\delta)\rbrace$,
\item[(3)] For any $G\in\mathcal{G}$, there exists $j\in \lbrace 1,\dots,  N_B(\delta)\rbrace$ such that $G_j\subset G\subset H_j$.
\end{itemize}

We begin with a general upper bound when $\mathcal{G}$ has a given entropy rate. It allows to deduce easily Theorem \ref{prop:aminima} in Section \ref{s:upperbound}. 
\begin{prop}
\label{prop:gen}
Suppose $\mathcal{G}$ contains $G^\star_K$  and satisfies, for some $0<\rho<1$:
\beqn
\label{entropycond}
\mathcal{H}(\mathcal{G},d_\Delta,\delta)\leq c\delta^{-\rho}.
\eeqn
Let $\hat G_{n}$ the set introduced in \eqref{generalsetestimation} where $\mathcal{G}$ satisfies \eqref{entropycond}. Suppose the noise assumption is satisfied and consider a kernel $\mathcal{K}_\eta$ defined as in \eqref{dk} satisfying the Kernel assumption. Then, we have:
$$
\E d_{f,g}(\hat G_{n,m},G_K^\star)\leq C\inf_{\lambda\in\R^d_+}\left[\left(\frac{\Pi_{i=1}^d\lambda_i^{-\beta_i}}{\sqrt{n}}\right)^{\frac{2(\alpha+1)}{\alpha+2+\rho\alpha}} +\sup_{G\in\mathcal{G}}\left(R_K-R_K^\lambda\right)( G)+\sum_{i=1}^d (n\lambda)^{-\frac{1}{1+\rho}}\right],
$$
where $C>0$ is a generic constant.
\end{prop} 
Such a result is an extension of Theorem \ref{prop:aminima} for general set $\mathcal{G}$ with given entropy conditions. The main ingredient of the proof is a generalization of Lemma 5.11 in \cite{vdg}. The proof is postponed to Section \ref{s:appendix}. \\

Such a generality allows to deal with various constraints on the problem. In this paper, we deal with assumptions on the smoothness of the Bayes classifier boundary. Alternative constraints could be investigated. By the way, it may be possible to consider plug-in type assumption as in \cite{AT} or \cite{les_pinkfloyd}, convex sets or finite Vapnik Chervonenkis classes as in \cite{nedelec} or \cite{kolt}.


\section{Proofs}
\label{s:proof}

In this section, with a slight abuse of notations, $C,c,c'>0$ denotes generic constants that may vary from line to line, and even in the same line. Given two real sequences $(a_n)_{n\in\mathbb{N}}$ and $(a_n)_{n\in\mathbb{N}}$, the notation $a_n\approx b_n$ (resp. $a_n\lesssim b_n$) means that there exists generic constants $C,c>0$ such that $ca_n\leq b_n\leq Ca_n$ (resp. $a_n\leq Cb_n$) for all $n\in \mathbb{N}$.

\label{section:proofs}
\subsection{Proof of Theorem \ref{thm:lbboundary}}
 The proof starts as in \cite{mammen} but then uses some arguments which are specific to the inverse problem literature (see for instance \cite{Butucea} or \cite{meister}). 

Let $\mathcal{F}_1$ a finite class of densities and $g_0$ a fixed density such that $(f,g_0)\in \mathcal{F}_{\mathrm{frag}}$ for all $f\in\mathcal{F}_1$. The contain of $\mathcal{F}_1$ and the value of $g_0$ will be precised later on. Then, for all estimator $\hat G_{n,m}$ of the set $G_K^{\star}$, we have
\begin{eqnarray}
\sup_{(f,g)\in \mathcal{F}_{frag}} \mathbb{E}_{f,g} &&\hspace{-0.9cm}d_{\Delta}(\hat G_{n,m},G_K^{\star})
\geq  \sup_{(f,g_0), f \in \mathcal{F}_1} \mathbb{E}_{f,g} d_{\Delta}(\hat G_{n,m},G_K^{\star}), \nonumber \\
& \geq & \mathbb{E}_{g_0} \left[ \frac{1}{\sharp\mathcal{F}_1} \sum_{f\in \mathcal{F}_1} \mathbb{E}_f \left\lbrace d_{\Delta}(\hat G_{n,m},G_K^{\star}) | X_1^{(2)}, \dots,X_m^{(2)} \right\rbrace \right]. 
\label{eq:borneinf}
\end{eqnarray}

\subsubsection{Construction of $\mathcal{F}_1$}
Concerning the density $g_0$, we deal with the uniform density on $[0,1]^2$, i.e.
$$ g_0(x) = \mathbf{1}_{\lbrace x \in [0,1]^2 \rbrace}, \forall x\in \R^2.$$
Now, we have to define the class $\mathcal{F}_1$. First, we consider a function $\varphi$ infinitely differentiable defined on $\R$ such that $\mathrm{supp}(\varphi)=[-1,1]$, $\varphi(t) \geq 0$ for all $t\in\R$ and $\| \varphi \|_{\infty} = \varphi(0)=1$. Let $M \geq 2$ an integer which will be allowed to depend on $n$ and $\tau>0$ a positive constant. Then, for all $j\in \lbrace 1,\dots,M \rbrace$, we set
$$ \varphi_j(t) = \tau M^{-\gamma} \varphi \left( M \left[ t - \frac{2j-1}{M} \right]\right), \ \forall t\in \R.$$
For all $\omega \in \lbrace 0,1 \rbrace^M$ and all $t\in \R$, we define
$$ b(t,\omega) = \frac{1}{2} + \sum_{j=1}^M \omega_j\varphi_j(t).$$
In the specific case where $\omega_j=1$ for all $j\in \lbrace 1,\dots,M \rbrace$, we write $b(t,\mathbf{1})$. Then, let $b_0$ and $C^{\star}$ positive constants which will be precised later on. We define the function $f_0: \R^2 \to \R$ as $f_0(x)=0$ for all $x\not \in [0,1]^2$ and
$$ f_0(x) = 
\left\lbrace
\begin{array}{l}
1 + 2\eta_0, \forall x_2 \in [0,1/2],\\
\\
1-\eta_0 - b_0, \forall x_2 \in [b(x_1,\mathbf{1}),1], \\
\\
1 + \left( \frac{b(x,\mathbf{1})-x_2}{c_2} \right)^{1/\alpha} - C^{\star}M^{-\gamma/\alpha}, \forall x_2 \in [1/2,b(x_1,\mathbf{1})],
\end{array}
\right.$$
where $ C^{\star} =3/2. (\tau/c_2)^{1/\alpha}$ and $b_0>0$ is such that $\int f_0(x)dx=3/4$. The condition on $C^\star$ ensures that $f_0(x)< 1$ for all $x_2 \in [1/2,b(x_1,\mathbf{1})]$. We will also use the function $f_1$ defined as 
$$ f_1(x) = 
\left\lbrace
\begin{array}{l}
0, \forall x\in [0,1]^2,\\
\\
\frac{b_1}{(1+x_2)^2. (1+x_1)^2} , \forall x\not\in [0,1]^2, 
\end{array}
\right.
$$
where $\mathcal{C}_1$ is such that $\int f_1(x)dx=1/4$. Finally, the set $\mathcal{F}_1$ will be defined as
$$ \mathcal{F}_1 = \left\lbrace f_{\omega}, \  \omega \in [0,1]^M  \right\rbrace,$$
where for a given $\omega \in \lbrace 0,1 \rbrace^M$,
\begin{equation}
f_{\omega}(x) = f_0(x) + f_1(x) + \sum_{j=1}^M \omega_j \rho_j(x).
\label{eq:fomega}
\end{equation}
for some functions $(\rho_j)_{j=1\dots M}$ which are explicited below. In order to complete the construction of the set $\mathcal{F}_1$, we have to provide a precise definition of the $\rho_j$ and to prove that the $f_{\omega}$ define probability density functions for all $\omega \in \lbrace 0,1\rbrace^M$. \\

We first start with the construction of the $\rho_j$. For all $x\in\R$, let $\rho: \R \to [0,1]$ the function defined as
$$ \rho(x) = \frac{1-\cos(x)}{\pi x^2}, \ \forall x\in \mathbb{R},$$
with associate Fourier transform $\mathcal{F}[\rho](t) = (1-|t|)_+$. In particular, $\mathrm{supp}\ \mathcal{F}[\rho] = [-1,1]$.  
For all $j\in \lbrace 1, \dots, M \rbrace$ and $x_2\in \R$, introduce
\begin{equation}
\rho_{(2)}(x_2) = \cos \left( \frac{x_2- 1/2(1+\tau M^{-\gamma})}{3/2\pi^{-1} \tau M^{-\gamma}} \right) \rho \left( \frac{x_2- 1/2(1+\tau M^{-\gamma})}{3\pi^{-1} \tau M^{-\gamma}} \right).
\label{eq:Gj2}
\end{equation}
By the same way, for all $j\in \lbrace 1, \dots, M \rbrace$, we define
\begin{equation}
\rho_{j,(1)}(x_1) = \cos \left[ \frac{\pi}{3}\left( \frac{x_1- j/M}{M^{-1}} \right)\right] \rho \left[  \frac{\pi}{6} \left(\frac{x_1- j/M}{M^{-1}} \right)\right].
\label{eq:Gj1}
\end{equation}
Then, for all $j\in \lbrace 1, \dots, M \rbrace$ and $x=(x_1,x_2) \in [0,1]^2$, we set 
\begin{equation}
\rho_j(x) = c^\star (\tau M^{-\gamma})^{1/\alpha} \ \rho_{(2)}(x_2) \rho_{j,(1)}(x_1),
\label{eq:Gj}
\end{equation}
for some constant $c^{\star}$ explicited below.  \\

Now, we prove that the $f_{\omega}$ introduced in (\ref{eq:fomega}) define density functions. First, remark that
$$ \sum_{j=1}^M | \rho_j(x) | \leq \left\lbrace
 \begin{array}{c}
  C M^{-\gamma/\alpha}(1+x_1)^{-2}(1+x_2)^{-2}, \ \forall x \not \in [0,1]^2,\\
  C M^{-\gamma/\alpha}  , \ \forall x \in [0,1]^2,
 \end{array}
\right.$$
This ensures that $f_{\omega}\geq 0$ for all $\omega\in \lbrace 0,1 \rbrace^M$, at least for $M$ large enough. Then recall that both $f_0$ and $f_1$ are designed in order to guarantee that $\int (f_0+f_1)(x)dx=1$. Hence, we only have to show that $\int \rho_j(x) dx =0$ for all $j\in \lbrace 1,\dots, M \rbrace$. In fact, it is only necessary to prove that $\int \rho_{(2)}(x_2) dx_2 =0$. First remark that $\int \rho_{(2)}(x_2) dx_2 = \int \tilde \rho_{(2)}(x_2) dx_2$ where $\tilde \rho_{(2)}(x_2) = \rho_{(2)}(x_2+1/2(1+\tau M^{-\gamma}))$ for all $x_2\in \R$. Then, using simple algebra
\beqnn \mathcal{F}[\rho_{(2)}](0) = \frac{1}{2} \mathcal{F}\left[\rho\left( \frac{.}{3 \pi^{-1} \tau M^{-\gamma}}\right)\right] \left( \pm \frac{1}{3/2 \pi^{-1}\tau M^{-\gamma}}  \right)\hspace{-0.4cm}&=&\hspace{-0.3cm}\frac{3}{2} \pi^{-1} \tau M^{-\gamma}\mathcal{F}[\rho] \left( \pm 2  \right)\\
&=&\hspace{-0.3cm} 0, 
\eeqnn
since the support of the Fourier transform of $\rho$ is $[-1;1]$. Hence, for all $\omega \in \lbrace 0,1 \rbrace^M$, $f_{\omega}$ is a density function. \\

In order to conclude the proof, we have to show that
\begin{equation} 
(f_{\omega},g_0) \in \mathcal{F}_{\mathrm{frag}} \ \forall \omega \in \lbrace 0,1 \rbrace^M,
\label{eq:asser1}
\end{equation}
which allows to use the bound (\ref{eq:borneinf}),
\begin{equation}
Q \left\lbrace x\in K: | f_{\omega}(x)-g_0(x) | \leq \eta \right\rbrace \leq c_2\eta^{\alpha} \ \forall \omega \in \lbrace 0,1 \rbrace^M \ \mathrm{and} \ \forall \eta\leq \eta_0, 
\label{eq:asser2}
\end{equation}
which means that the \textit{Margin assumption} is satisfied for our test functions and that
\begin{equation}
\mathbb{E}_{g_0} \mathbb{E}_{f_{\omega}} \left\lbrace d_{\Delta}(\hat G_{n,m},G_K^{\star}) | X_1^{(2)},\dots,X_m^{(2)} \right\rbrace  \geq C n^{-\frac{\gamma}{\gamma \left(\frac{2}{\alpha}+1\right) + 2\beta_1 + 2\beta_2 \gamma + 1}},
\label{eq:asser3}
\end{equation}
for some positive constant $C$.\\

\subsubsection{Main assumptions check}
We first start with the proof of (\ref{eq:asser1}). First remark that for all $j\in \lbrace 1, \dots, M \rbrace$, the function $\rho_{j}(.)$ is bounded from above by $C M^{-\gamma/\alpha}$ for some $C>0$. Then, using simple algebra
\begin{eqnarray*}
x_2 \in [1/2; b(x_1,\mathbf{1})]
& \Rightarrow & \frac{1}{2} \leq x_2 \leq \frac{1}{2} + \tau M^{-\gamma},\\
& \Rightarrow & -\frac{\tau M^{-\gamma}}{2} \leq x_2 - \frac{1}{2} - \frac{\tau M^{-\gamma}}{2} \leq  \frac{\tau M^{-\gamma}}{2},\\
& \Rightarrow & - \frac{\pi}{6} \leq \frac{x_2- 1/2(1+\tau M^{-\gamma})}{3\pi^{-1} \tau M^{-\gamma}} \leq \frac{\pi}{6},\\
& \Rightarrow & \rho_{(2)} (x_2) \geq \frac{9}{4\pi^3}.
\end{eqnarray*}
The same kind on minoration holds for the function $\rho_{j,(1)}$. Hence the $\rho_j$ are uniformly bounded from below on $[1/2;b(x_1,\mathbf{1})]$. For all $\omega \in \lbrace 0,1 \rbrace^M$ and for all $x \in [0,1]^2$, we then have
$$ f_{\omega}(x) \geq 1 + \left( \frac{b(x,\mathbf{1})-x_2}{c_2} \right)^{1/\alpha} \geq g_0(x), \ \forall x_2 \in [1/2,b(x_1,\omega)],$$
for $c^\star$ large enough. This ensures that
$$ \lbrace x \in [0,1]^2: f_{\omega}(x) \geq g_0(x) \rbrace = \lbrace x \in [0,1]^2: 0 \leq x_2 \leq b(x_1, \omega) \rbrace .$$
In order to conclude the proof of (\ref{eq:asser1}), we only have to remark that the function $b(., \omega)$ belongs to $\Sigma(\gamma,L)$ for all $\omega \in \lbrace 0,1 \rbrace^M$, at least for $M$ small enough.\\

Now, we consider the margin assumption (\ref{eq:asser2}). First, we consider the case where $\eta < [\tau c_2^{-1}]^{1/\alpha} M^{-\gamma/\alpha}<\eta_0$. Clearly, following our choices of $b_0$ and $C^{\star}$, we have that 
$$ | f_{\omega}(x)-g_0(x) | \leq \eta \Rightarrow x_2 \in [1/2; b(x_1,\omega)] \Rightarrow x_2 \leq b(x_1,\omega).$$
Moreover, for all $x\in [0,1]^2$ such that $x_2 \leq b(x_1,\omega)$, we have
$$ (f_{\omega} - g_0)(x) =  \left( \frac{b(x,\mathbf{1})-x_2}{c_2} \right)^{1/\alpha} + \sum_{j=1}^M \omega_j \rho_j(x) - C^{\star}M^{-\gamma/\alpha},$$
where                      
$$ \sum_{j=1}^M \omega_j \rho_j(x) - C^{\star} M^{-\gamma/\alpha} > 0, \ \forall x_2 \in \left[ \frac{1}{2},b(x_1,\omega) \right].$$
Thus
$$ | f_{\omega}(x)-g_0(x) | \leq \eta \Rightarrow \left( \frac{b(x,\omega)-x_2}{c_2} \right)^{1/\alpha} \leq \eta \Rightarrow x_2 \geq b(x_1,\omega) - c_2 \eta^{\alpha},$$
which proves the margin assumption when $\eta < [\tau c_2^{-1}]^{1/\alpha} M^{-\gamma/\alpha}$. Now, in the case where $\eta_0 >\eta > [\tau c_2^{-1}]^{1/\alpha} M^{-\gamma/\alpha}$, we have
$$ | f_{\omega}(x)-g_0(x) | \leq \eta \Rightarrow 1/2 < x_2 <b(x_1, \mathbf{1}),$$
which entails
$$ Q \left\lbrace x\in K: | f_{\omega}(x)-g_0(x) | \leq \eta \right\rbrace  \leq \tau M^{-\gamma} \leq c_2 \eta^{\alpha}.$$
This concludes this part.\\

\subsubsection{Final minoration}
Now, we can deal with the lower bound (\ref{eq:asser3}). The proof is based on classical tools which can be found for instance in \cite{booktsybakov}, \cite{mammen}, \cite{Butucea} or \cite{meister}. First remark that the shape of $G_K^{\star}$ depends on the value of $\omega$. For the sake of convenience, we omit the dependency with respect to this quantity. For all $\omega\in \lbrace 0,1 \rbrace ^M$, recall that
$$ G_K^{\star} = \lbrace x\in [0,1]^2: f_{\omega}(x) \geq g_0(x) \rbrace = \lbrace x\in [0,1]^2: 0 \leq x_2 \leq b(x_1,\omega) \rbrace.$$
Using Assouad Lemma and classical tools designed for instance in \cite{booktsybakov}, we get
\begin{equation}
\mathbb{E} \left[ d_{\Delta}(\hat G_{n,m}, G_K^{\star}) | Y_1,\dots,Y_m \right] \geq \frac{M}{2} \| \varphi_1 \|_1 \int \min \left[ dP_{11},dP_{10} \right],
\label{eq:borninfinter}
\end{equation}
where $P_{11}$ denotes the law of $(Z_i^{(1)})_{i=1\dots n}$ when the density of the $X_i^{(1)}$ is $f_{\omega_{11}}$. In the following, we will choose $M$ in order to guarantee that the term $\int \min \left[ dP_{11},dP_{10} \right]$ is bounded from below. Consequently, the lower bound will be determined by the corresponding value of $M \| \varphi_1 \|_1$. Since the observations are independent
$$\int \min \left[ dP_{11},dP_{10} \right] \geq 1 - \sqrt{ \left( 1+ \chi^2(P_1,P_0) \right)^n -1},$$
where $\chi^2(P_a,P_b)$ denotes the chi-square divergence between two given probability measures $P_a$ and $P_b$, and $P_0,P_1$ are the law of the variable $Z_1^{(1)}=X_1^{(1)}+\epsilon_1^{(1)}$ when the density of the $X_i^{(1)}$ is respectively $f_{\omega_{11}}$ or $f_{\omega_{10}}$. In the following, our aim is to find a satisfying upper bound for $\chi^2(P_1,P_0)$.\\

First, remark that we can find $\tilde c>0$ such that for all $x \not \in [0,1]^2$ and all $\omega\in \lbrace 0,1 \rbrace^M$, $ f_{\omega}(x) \geq \tilde c f_1(x)$. Hence, using simple algebra, we get that
\begin{equation}
f_{\omega}*\eta(x) \geq \frac{C}{(1+x_1^2)(1+x_2^2)}, \ \forall x\in \mathbb{R}^2,
\label{eq:mino}
\end{equation}
for some $C>0$. In the following, given $f,\eta_1$ and $\eta_2$, we denote by $f*\eta$ the convolution product in dimension two, i.e.
$$ f*\eta(x) = \int_{\mathbb{R}}\int_{\mathbb{R}} f(x_1-y_1,x_2-y_2)\eta_1(y_1)\eta_2(y_2)dy_1dy_2, \ \forall x\in \mathbb{R}^2.$$ 
Then, using (\ref{eq:fomega}) and (\ref{eq:mino}),
\begin{eqnarray*}
\chi^2(P_1,P_0)
& = & \int_{\mathbb{R}}\int_{\mathbb{R}} \frac{\lbrace (f_{\omega_{11}}-f_{\omega_{10}})*\eta(x) \rbrace^2}{f_{\omega_{11}}*\eta(x)} dx,\\
& \leq & C \int_{\mathbb{R}}\int_{\mathbb{R}} (1+x_1^2) (1+x_2^2) \lbrace \rho_1*\eta(x) \rbrace^2 dx. 
\end{eqnarray*}
Hence
\begin{eqnarray*}
\chi^2(P_1,P_0)
& \leq & C \int_{\mathbb{R}}\int_{\mathbb{R}} \lbrace \rho_1*\eta(x) \rbrace^2 dx + C \int_{\mathbb{R}}\int_{\mathbb{R}} x_2^2 \lbrace \rho_1*\eta(x) \rbrace^2 dx \\
& & + C \int_{\mathbb{R}}\int_{\mathbb{R}} x_1^2 \lbrace \rho_1*\eta(x) \rbrace^2 dx + C \int_{\mathbb{R}}\int_{\mathbb{R}} x_1^2 x_2^2 \lbrace \rho_1*\eta(x) \rbrace^2 dx, \\
& := & A_1 + A_2 + A_3 + A_4,
\end{eqnarray*}
where the $\rho_j$ are defined in (\ref{eq:Gj}). In the following, we only consider the bound of $A_1$, the other terms being controlled in the same way. We get
\begin{eqnarray*}
A_1 & = & C \int_{\mathbb{R}}\int_{\mathbb{R}} \lbrace \rho_1*\eta(x) \rbrace^2 dx,\\
& = & C M^{-2\gamma/\alpha} \int_{\mathbb{R}}\int_{\mathbb{R}} \left\lbrace \int_{\mathbb{R}}\int_{\mathbb{R}} \rho_{(2)}(x_2-y_2) \rho_{j,(1)}(x_1-y_1) \eta_1(y_1) \eta_2(y_2)dy_1 dy_2 \right\rbrace^2 dx,\\
& = & C M^{-2\gamma/\alpha} \int_{\mathbb{R}}\int_{\mathbb{R}} | \mathcal{F}[\rho_{(2)}](t_2) |^2 | \mathcal{F}[\rho_{1,(1)}](t_1) |^2 |\mathcal{F}[\eta_1](t_1)|^2 |\mathcal{F}[\eta_2](t_2)|^2 dt_1dt_2,\\
& = & C M^{-2\gamma/\alpha} A_{1,1} A_{1,2},
\end{eqnarray*}
where
$$ A_{1,1} = \int_{\mathbb{R}} | \mathcal{F}[\rho_{(1)}](t_1) |^2 |\mathcal{F}[\eta_1](t_1)|^2 dt_1, \ A_{1,2}=\int_{\mathbb{R}}\int_{\mathbb{R}} | \mathcal{F}[\rho_{1,(2)]}(t_2) |^2  |\mathcal{F}[\eta_2](t_2)|^2 dt_2,$$
and $\rho_{(1)}$, $\rho_{1,(2)}$ are respectively defined in (\ref{eq:Gj2}),(\ref{eq:Gj1}). We first deal with the term $A_{1,2}$. Using simple algebra, we get
\begin{eqnarray*}
A_{1,2}
& = & \int_{\mathbb{R}} | \mathcal{F}[\rho_1^{(1)}](t_1) |^2 |\mathcal{F}[\eta_1](t_1)|^2 dt_1,\\
& = & \int_{\mathbb{R}} \left| \mathcal{F}\left[\rho\left( \frac{.}{3\pi^{-1} \tau M^{-\gamma}} \right)\right] \left(t_1 \pm \frac{1}{3/2 \pi^{-1} \tau M^{-\gamma}}  \right) \right|^2 |\mathcal{F}[\eta_1](t_1)|^2 dt_1,\\
& = & (3\pi^{-1})^2 \tau^2 M^{-2\gamma} \int_{\mathbb{R}} \left| \mathcal{F}[\rho] \left(3\pi^{-1}\tau M^{-\gamma} t_1 \pm \frac{3}{3/2}  \right) \right|^2 |\mathcal{F}[\eta_1](t_1)|^2 dt_1. 
\end{eqnarray*}
Then, setting $s_1= 3\pi^{-1} \tau M^{-\gamma} t_1$ and using the \textit{Noise assumption}, we obtain
\begin{eqnarray*}
A_{1,2} 
& = & 3\pi^{-1} \tau M^{-\gamma} \int_{\mathbb{R}} | \mathcal{F}[\rho] (s_1 \pm 2) |^2 \left| \mathcal{F}[\eta_1]\left( \frac{s_1}{3\pi^{-1}\tau M^{-\gamma}} \right) \right|^2 ds_1, \\
& = & 3\pi^{-1} \tau M^{-\gamma} \int_{1}^3 | \mathcal{F}[\rho] (s_1 \pm 2) |^2 \left| \mathcal{F}[\eta_1]\left( \frac{s_1}{3\pi^{-1}\tau M^{-\gamma}} \right) \right|^2 ds_1, \\
& \leq & C M^{-\gamma-2\beta_2\gamma} \int_1^3 | \mathcal{F}[\rho] (s_1 \pm 2) |^2 |s_1|^{-2\beta_1} ds_1,\\
& \leq & C M^{-\gamma-2\beta_2 \gamma}.
\end{eqnarray*}
Using a similar algebra for the term $A_{1,1}$, we obtain
$$ A_{1,2} \leq C  M^{-1-2\beta_1}.$$
Similar bounds are available for $A_2,A_3$ and $A_4$ since $\mathcal{F}[\rho]$ and its weak derivative are bounded by $1$ and supported on $[-1;1]$. In particular, we use the fact that for all $t\in \R$
$$ \mathcal{F}[\rho_{1,(2)}](t) = 3 \pi^{-1} \tau M^{-\gamma} \mathcal{F}[\rho](3\pi^{-1} \tau M^{-\gamma} t \pm 2 ),$$
and 
$$ \frac{d}{dt} \mathcal{F}[\rho_{1,(2)}](t) = - i (3 \pi^{-1} \tau M^{-\gamma})^2 t. \mathcal{F}[\rho](3\pi^{-1} \tau M^{-\gamma} t \pm 2 ),$$ 
for all $t$ in a subset of $\R$ having a Lebesgue measure equal to $1$. \\

The above equations lead to the following upper bound:
$$ \chi^2(P_1,P_0) \leq C M^{-\gamma (2/\alpha + 1) - 2\beta_1 \gamma - 2\beta_2 -1}.$$
Then, $\chi^2(P_1,P_0) \leq C/n$ for some constant $C>0$ as soon as
$$ M= M_n \sim n^{\frac{1}{\gamma (2/\alpha + 1) + 2\beta_1 + 2\beta_2\gamma +1}}.$$
Finally, going back to equation (\ref{eq:borninfinter}), we obtain
\begin{eqnarray*}
\mathbb{E} \left[ d_{\Delta}(\hat G_{n,m}, G_K^{\star}) | Y_1,\dots,Y_m \right] 
& \geq & \frac{M_n}{2} \| \varphi_1 \|_1 \int \min \left[ dP_{11},dP_{10} \right],\\
& \geq & C M_n \| \varphi_1 \|_1,\\
& = & C\tau M_n^{-\gamma} \int_0^1 \varphi_1(t)dt,\\
& \sim &  M_n^{-\gamma} = n^{-\frac{\gamma}{\gamma (2/\alpha + 1) + 2\beta_1  + 2\beta_2\gamma +1}},
\end{eqnarray*}
which concludes the proof.
\begin{flushright}
$\Box$
\end{flushright}

\subsection{Proof of Proposition \ref{prop:gen}}
Consider the empirical processes $\nu_n^{(j)}$, for $j\in\{1,2\}$, defined as:
\beqn
\label{epj}
\nu_n^{(j)}(G)=\frac{1}{\sqrt{n}}\sum_{i=1}^n \left[ h_{G,\lambda}(Z_i^{(j)})-\E h_{G,\lambda}(Z^{(j)})\right],
\eeqn
where the $h_{G,\lambda}(.)$ have been introduced in (\ref{hG}). In particular, remark that  for all $i\in \lbrace 1,\dots, n \rbrace$, $G \subset K$, 
\begin{eqnarray*} 
\mathbb{E}\left[ h_{G,\lambda}(Z_i^{(1)}) \right] 
& = & \int_G \frac{1}{\lambda} \mathbb{E}\left[ \mathcal{K}_\eta \left( \frac{X_i^{(1)}+\epsilon_i^{(1)} - x}{\lambda} \right) \right] dx,\\
& = & \int_G \frac{1}{\lambda} \mathbb{E}\left[ \mathcal{K} \left( \frac{X_i^{(1)} - x}{\lambda} \right) \right] dx = \int_{\mathbb{R}^d} f(x) \mathcal{K}_\lambda*\mathbf{1}_{\lbrace . \in G\rbrace }(x)dx.  
\end{eqnarray*}
Hence, using (\ref{derm2}), we can write
\begin{eqnarray}
\lefteqn{\int (f-g)(\mathcal{K}_\lambda*\mathbf{1}_{\lbrace . \in G^\star_K \rbrace }- \mathcal{K}_\lambda*\mathbf{1}_{\lbrace . \in \hat G_n \rbrace })} \nonumber \\
&\leq &\frac{1}{\sqrt{n}}(\nu_n^{(1)}(G^{\star C})-\nu_n^{(1)}(\hat{G}^{\lambda C}_{n}))+ \frac{1}{\sqrt{n}}(\nu_n^{(2)}(G^{\star})-\nu_n^{(2)}(\hat{G}^{\lambda}_{n}).
\label{eq:prop1_inter}
\end{eqnarray}
Now denoting $\Lambda=\Pi_{i=1}^d\lambda_i^{-\beta_i-\frac{1}{2}}$, $c(\lambda)=\Pi_{i=1}^d\lambda_i^{-\beta_i}$ and $\rho=2/\gamma$, consider the event
\beqnn
\Omega := \{d_{\Delta}(\hat{G}_n,G^\star_K)\geq c(\lambda)^{-\frac{2}{1+\rho}}n^{-\frac{1}{1+\rho}}\Lambda^{\frac{2}{1+\rho}}\}.
\eeqnn
If the event $\Omega$ holds, using Lemma 2 of \cite{mammen}, we get
\begin{eqnarray*}
\lefteqn{\int (f-g)(\mathcal{K}_\lambda*\mathbf{1}_{\lbrace . \in G^\star_K \rbrace }- \mathcal{K}_\lambda*\mathbf{1}_{\lbrace . \in \hat G_n \rbrace })}\\
&\leq & \frac{d_{\Delta}^{\frac{1-\rho}{2}}(\hat{G}^\lambda_{n,m},G^\star)c(\lambda)}{\sqrt{n}}\left[\frac{\nu^{(1)}_n(G^{\star C})-\nu^{(1)}_n(\hat{G}^{\lambda C}_{n,m})}{c(\lambda)d_{\Delta}^{\frac{1-\rho}{2}}(\hat{G}^\lambda_{n,m},G^\star)\vee c(\lambda)^{\frac{2\rho}{(1+\rho)}}n^{-\frac{1-\rho}{2+2\rho}}\Lambda^{\frac{1-\rho}{1+\rho}}}\right.\\
&& \hspace{1.5cm}   +\left.\frac{\nu^{(2)}_n(G^{\star})-\nu^{(2)}_n(\hat{G}^{\lambda}_{n,m})}{c(\lambda)d_{\Delta}^{\frac{1-\rho}{2}}(\hat{G}^\lambda_{n,m},G^\star)\vee c(\lambda)^{\frac{2\rho}{(1+\rho)}}n^{-\frac{1-\rho}{2+2\rho}}\Lambda^{\frac{1-\rho}{1+\rho}}}\right],\\
& \leq & \frac{d_{f,g}^{\frac{1-\rho}{2}\frac{\alpha}{\alpha+1}}(\hat{G}^\lambda_{n,m},G^\star)c(\lambda)}{\sqrt{n}}[V^{(1)}_n+V^{(2)}_n],\\
\end{eqnarray*}
where for $j\in\{1,2\}$, $V^{(j)}_n$ is the random variable defined as
\beqn
\label{vnj}
V^{(j)}_n=\sup_{G\in\GG}\frac{|\nu^{(j)}_n(G^\star)-\nu^{(j)}_n(G)|}{c(\lambda)\no \mathbf{1}_G-\mathbf{1}_{G^*}\no_{2,X^{(j)}}^{1-\rho}\vee c(\lambda)^{\frac{2\rho}{(1+\rho)}}n^{-\frac{1-\rho}{2+2\rho}}\Lambda^{\frac{1-\rho}{1+\rho}}}.
\eeqn
Lemma \ref{noisyep} in Section \ref{s:preproof} shows that the variable $V_n^{(1)}+V_n^{(2)}$ has controlled moments. Indeed, the bracketing entropy related to the set $\mathcal{G}$ is $\rho = (d-1)/\gamma=1/\gamma$. Using the Young's inequality $xy^r\leq ry+(1-r)x^{1/(1-r)}$ with $r=\frac{1-\rho}{2}\frac{\alpha}{\alpha+1}$, we get
\begin{eqnarray}
\lefteqn{\int (f-g)(\mathcal{K}_\lambda*\mathbf{1}_{\lbrace . \in \hat G_n \rbrace }-\mathcal{K}_\lambda*\mathbf{1}_{\lbrace . \in G^\star_K \rbrace })} \nonumber \\
& \leq & c\left(\frac{c(\lambda)}{\tau^{-1}\sqrt{n}}[V_n^{(1)}+V_n^{(2)}]\right)^{\frac{2(\alpha+1)}{\alpha+2+\rho\alpha}} + \tau d_{f,g}(\hat G,G_K^\star).
\label{varcontrolonomega}
\end{eqnarray}
Note that
\beqnn
d_{f,g}(\hat G,G^\star)&=&\int (f-g)(\mathcal{K}_\lambda*\mathbf{1}_{\lbrace . \in \hat G_n \rbrace }-\mathcal{K}_\lambda*\mathbf{1}_{\lbrace . \in G^\star_K \rbrace })+\left(R_K-R_K^\lambda\right)(\hat G,G^\star)\\
&\leq & \int (f-g)(\mathcal{K}_\lambda*\mathbf{1}_{\lbrace . \in \hat G_n \rbrace }-\mathcal{K}_\lambda*\mathbf{1}_{\lbrace . \in G^\star_K \rbrace })+2\sup_{G\in\mathcal{G}}\left(R_K-R_K^\lambda\right)(\hat G).
\eeqnn
From above, we have coarselly:
\begin{eqnarray*}
d_{f,g}(\hat G,G^\star)\mathbf{1}_\Omega\leq \left(\frac{1}{1-\tau}\right)c\left(\frac{c(\lambda)}{\tau^{-1}\sqrt{n}}[V_n^{(1)}+V_n^{(2)}]\right)^{\frac{2(\alpha+1)}{\alpha+2+\rho\alpha}} +2\sup_{G\in\mathcal{G}}\left(R_K-R_K^\lambda\right)(G).
\end{eqnarray*}
In order to end up the proof, let us consider the following decomposition:
\begin{eqnarray}
d_{f,g}(\hat G,G^\star) & = & d_{f,g}(\hat G,G^\star)\mathbf{1}_{\Omega}+d_{f,g}(\hat G,G^\star)\mathbf{1}_{\Omega^C}
\label{decfinal}
\end{eqnarray}
Moreover, note that on the event $\Omega^C$, we have:
$$
d_{\Delta}(\hat{G}_n,G^\star_K)\leq c(\lambda)^{-\frac{2}{1+\rho}}n^{-\frac{1}{1+\rho}}\Lambda^{\frac{2}{1+\rho}}=(n\lambda)^{-\frac{1}{1+\rho}}.
$$
Hence, we can conclude that
\begin{eqnarray}
d_{f,g}(\hat G,G^\star)\leq c_1\left(\frac{c(\lambda)}{\sqrt{n}}[V_n^{(1)}+V_n^{(2)}]\right)^{\frac{2(\alpha+1)}{\alpha+2+\rho\alpha}} +c_2\sup_{G\in\mathcal{G}}\left| R_K-R_K^\lambda\right| ( G)+c_3(n\lambda)^{-\frac{1}{1+\rho}}.
\label{fin_prop1}
\end{eqnarray}
Integrating the last inequality, we get the result of Theorem \ref{prop:aminima}.
\begin{flushright}
$\Box$
\end{flushright}

\subsection{Proof of Theorem \ref{prop:aminima}}
The proof is a direct consequence of Proposition \ref{prop:gen} as soon as we remark that
$$ \mathcal{H} (\mathcal{G}(\gamma,L),d_\Delta,\delta) \leq c\delta^{-\frac{d-1}{\gamma}}, \ \forall \ \delta>0.$$
This result can be found in \cite{wvdv}.

\subsection{Proof of Corollary \ref{cor:rough}}
Thanks to the previous proof, we only have to propose a bound for the term 
$$\sup_{G\in\mathcal{G}}\left| R_K-R_K^\lambda\right| ( G).$$
Let $G\in \mathcal{G}$ be fixed. For the sake of convenience, we restrict ourselves to the particular case where $d=2$. The generalization to larger dimension is straightforward. Moreover, we restrict ourselves to the control of the bias term over  the compact $K'=[\epsilon,1-\epsilon]^{d-1}\times [0,1]$, where $\epsilon>0$ is a small positive constant chosen later on to have:
\beqn
\left|(R_{K}-R_{K}^\lambda)(G)-(R_{K'}-R_{K'}^\lambda)(G)\right|\leq C \psi_n(\alpha,\gamma,\beta),
\label{eq:rigbias}
\eeqn
where $\psi_n(\alpha,\gamma,\beta)$ is the expected rate of convergence. Using \eqref{eq:rigbias}, we can conduct the proof of Corollary \ref{cor:rough} over $K'$. Then, using Fubini, remark that, provided that $x_1-\lambda_1z_1\in[0,1]$:
\begin{eqnarray*}
(R_{K'}-R_{K'}^\lambda) ( G)
& = & \int_{K'} (f-g) \left( \mathcal{K}_\lambda * \mathbf{1}_G - \mathbf{1}_G \right) d\lambda\\
& = & \int_{K'} (f-g)(x) \left( \int_{\mathbb{R}^2} \frac{1}{\lambda} \mathcal{K} \left( \frac{x-z}{\lambda}\right) \left[ \mathbf{1}_G(z) - \mathbf{1}_G(x) \right] dz \right)dx\\
& = & \int_{K'} (f-g)(x) \left( \int_{\mathbb{R}^2}  \mathcal{K} \left( z \right) \left[ \mathbf{1}_G(x-\lambda z) - \mathbf{1}_G(x) \right] dz \right)dx,\\
& = &   \int_{\mathbb{R}^2}  \mathcal{K}\left( z \right) \int_{K'} (f-g)(x)  \left[ \mathbf{1}_G(x-\lambda z) - \mathbf{1}_G(x) \right]  dx dz.
\end{eqnarray*}
Note that if $K'=[\epsilon,1-\epsilon]\times [0,1]$ and $\mathrm{supp} \mathcal{K}_1=[-M,M]$, for any $\lambda_1\leq \lambda_\mathrm{max}$, the choice of $\epsilon=M\lambda_{\mathrm{max}}$ ensures that:
$$
x\in K'\Rightarrow 0\leq x_1-\lambda_1z_1\leq 1.
$$ 
Moreover, since
$$ x\in G \Leftrightarrow 0 \leq x_2 \leq b(x_1),$$
we get:
$$ x-\lambda z \in G \Leftrightarrow 0 \leq x_2 \leq\min(1, b(x_1 - \lambda_1 z_1)+ \lambda_2 z_2 ).$$
Finally, for all $z_1 \in \mathbb{R}$, since $b\in\Sigma(\gamma,L)$:
\begin{equation}
b(x_1 - \lambda_1 z_1) =p_{b,x_1}(x_1-\lambda_1 z_1) + \mathcal{O}(|\lambda_1z_1|^\gamma).
\label{eq:dvp_taylor}
\end{equation}
Hence, we obtain, using the crude bound $\no f-g\no_\infty\leq 2c_1$ and the assumptions over the kernel $\mathcal{K}$:
\begin{eqnarray}
\lefteqn{(R_{K'}-R_{K'}^\lambda) ( G)} \nonumber \\
& = &   \int_{\mathbb{R}^2} \hspace{-0.2cm} \mathcal{K}(z)\hspace{-0.2cm} \int_\epsilon^{1-\epsilon} \left[ \int_{\lambda_2z_2}^{\min(1,b(x_1-\lambda_1z_1)+\lambda_2z_2)} (f-g)(x)dx_2-\int_{0}^{b(x_1)} (f-g)(x)dx_2 \right]dx_1  dz \label{eq:biais_inter} \\
& \leq & 2c_1 \int_{\mathbb{R}^2}  \left| \mathcal{K}(z) \right| \int_\epsilon^{1-\epsilon} \left| p_{b,x_1}(x_1-\lambda_1z_1)-b(x_1)+\mathcal{O}(|\lambda_1z_1|^\gamma) + 2\lambda_2z_2\right| dx_1 dz \nonumber \\
& \leq & C(\lambda_1 + \lambda_2),
\label{biais_pourri}
\end{eqnarray}
where $C>0$ is a generic constant. Using (\ref{fin_prop1}) and (\ref{biais_pourri}), we obtain
$$ d_{f,g}(\hat G,G^\star)\leq c_1\left(\frac{c(\lambda)}{\sqrt{n}}[V_n^{(1)}+V_n^{(2)}]\right)^{\frac{2(\alpha+1)}{\alpha+2+\rho\alpha}}+ C(\lambda_1 + \lambda_2)+c_3(n\lambda)^{-\frac{1}{1+\rho}}.$$
We can conclude the proof with an appropriate choice for $\lambda_1$ and $\lambda_2$, noting that, in dimension $d=2$ for simplicity, \eqref{eq:rigbias} holds since for any $\lambda_1\leq\lambda_{\mathrm{max}}$:
\beqnn
&&\left|(R_{K}-R_{K}^\lambda)(G)-(R_{K'}-R_{K'}^\lambda)(G)\right|\leq \left|\int_{K\setminus K'}(f-g)\left(\mathbf{1}_G-\mathcal{K}_\lambda*\mathbf{1}_G\right)d\lambda\right|\nonumber\\
&\hspace{-0.5cm}= &\hspace{-0.5cm}  \left|\int_{\R^2} \mathcal{K}(z)\left(\int_0^{\epsilon}+\int_{1-\epsilon}^{1}\right)\left[ \int_{\lambda_2z_2}^{b(x_1-\lambda_1z_1)+\lambda_2z_2} (f-g)(x)dx_2-\int_{0}^{b(x_1)} (f-g)(x)dx_2 \right]dx_1  dz\right|\nonumber\\
&\leq &2\epsilon \leq C (\lambda_1 + \lambda_2)\nonumber,
\eeqnn
for $\epsilon=\lambda_{\mathrm{max}}M$.
\begin{flushright}
$\Box$
\end{flushright}

%
%
%
%
%
%
%
%
%
%
%
%
\section{Appendix}
\label{s:appendix}
\subsection{Technical Lemmas}

\begin{lemma}
\label{lip}
Let $Z$ be a random variable having density $f*\eta$ w.r.t. the Lebesgue measure. Assume that $\eta$ satisfies the \textit{Noise assumption} and that \textbf{(K1)} and \textbf{(K2)} hold. Then we have, 
\begin{eqnarray*}
(i) & & \E [h_{G,\lambda}(Z)-h_{G',\lambda}(Z)]^2\leq C d_\Delta(G,G') \prod_{i=1}^d \lambda_i^{-2\beta_i}.\\
(ii) & & \sup_{x\in K} | h_{G,\lambda}(x) - h_{G',\lambda}(x) | \leq C \prod_{i=1}^d \lambda_i^{-\beta_i-1/2},
\end{eqnarray*}
where $C>0$ is a generic constant.
\end{lemma}
\textsc{proof}
For the sake of convenience, we only consider the case where $d=1$. We first prove $(i)$. We have
\begin{eqnarray*}
 &&\E [h_{G,\lambda}(Z)-h_{G',\lambda}(Z)]^2
 = \\
&&\int_K \left[ \int_{\mathbb{R}} \frac{1}{\lambda} \mathcal{K}_\eta \left(\frac{z-x}{\lambda}\right) (\mathbf{1}_{\lbrace x\in G\rbrace}-\mathbf{1}_{\lbrace x\in G'\rbrace})\mathbf{1}_{\lbrace x\in K\rbrace} dQ(x) \right]^2 f*\eta(z)dz, \\
& \leq & c\int_{\mathbb{R}} \frac{1}{\lambda^2} \left| \mathcal{F}[\mathcal{K}_\eta(./\lambda)](t) \right|^2 \left| \mathcal{F}[ (\mathbf{1}_{\lbrace . \in G\rbrace}-\mathbf{1}_{\lbrace . \in  G'\rbrace})\mathbf{1}_{\lbrace . \in K\rbrace}](t) \right|^2 dt,\\
& \leq & C\lambda^{-2\beta} \int_{K} \mathbf{1}_{\lbrace t \in G\Delta G'\rbrace} dt,\\
& \leq & C\lambda^{-2\beta} d_{\Delta}(G,G').
\end{eqnarray*}
Indeed, for all $s\in\mathbb{R}$, using \textbf{(K3)}:
\begin{equation}
\frac{1}{\lambda^2} \left| \mathcal{F}[\mathcal{K}_\eta(./\lambda)](s) \right|^2 =  \left| \mathcal{F}[\mathcal{K}_\eta](s\lambda) \right|^2  \leq \sup_{t\in \mathbb{R}} \left| \mathcal{F}[\mathcal{K}_\eta](t) \right|^2 \leq C \lambda^{-2\beta},
\label{eq:inter}
\end{equation}
By the same way,
\begin{eqnarray*}
\sup_{x\in \mathbb{R}} | h_{G,\lambda}(x) - h_{G',\lambda}(x) |
& = & \sup_{x\in\mathbb{R}}  \int_{G\Delta G'} \frac{1}{\lambda} \left| \mathcal{K}_\eta \left(\frac{z-x}{\lambda} \right)\right| dx ,\\
& \leq & \sup_{x\in\mathbb{R}}  \int_{K} \frac{1}{\lambda} \left| \mathcal{K}_\eta \left(\frac{z-x}{\lambda} \right)\right| dx , \\
& \leq & C \sup_{x\in\mathbb{R}} \sqrt{ \int \frac{1}{\lambda^2} \mathcal{K}^2_\eta \left(\frac{z-x}{\lambda} \right) dx} \leq \lambda^{-\beta-1/2},
\end{eqnarray*}
where the last line is inspired by (\ref{eq:inter}).
\begin{flushright}
$\Box$
\end{flushright}


\subsection{Noisy Empirical process theory}
In this paragraph, we present the main ingredient for the proof of Theorem \ref{prop:aminima} and Proposition \ref{prop:gen}. We intend to analyze the behaviour of the increments of a noisy empirical process related with error measurements. The framework is much more general than model  \eqref{eq:observation} and deconvolution classifier of Section \ref{s:model}. Let us fix some notations. \\
Given a class of functions $\mathcal{G}$, we study the following risk minimization problem:
$$
g^\star=\arg\min_{g\in\mathcal{G} }R(g),
$$
where we have at our disposal a training sample $Z_1,\ldots ,Z_n$ of i.i.d. random variable with law $P_Z$. It differs from the law of $X$, denoted by $P_X$ and since $R(g):=\E[g(X)]$, we are faced to an inverse problem. For this purpose, we consider a indirect ERM procedure that can be written:
$$
\hat{g}_n^\lambda:=\arg\min_{g\in\mathcal{G}}\frac{1}{n}\sum_{i=1}^ng_\lambda(Z_i),
$$
where $g_\lambda:=\Phi_\lambda(g)$ is a smoothed version of $g\in\mathcal{G}$ and $\lambda\in\Lambda$ is a smoothing parameter (see the particular case $\Phi_\lambda(G)=h_G^\lambda$ in Section \ref{s:model}). 
 
We are interested in the control of the variance of $\hat{g}_n^\lambda$, which is equivalent to the study of the increments of the empirical process $\nu_n$ defined as:
\begin{equation}
\nu_n^\lambda(g)=\frac{1}{\sqrt{n}}\sum_{i=1}^n [g_\lambda(Z_i)-\E g_\lambda(Z)].
\label{eq:proc}
\end{equation}
The empirical process \eqref{eq:proc} is indexed by the set $\mathcal{G}^\lambda=\{g_\lambda=\Phi_\lambda(g),\,g\in\mathcal{G}\}$. We are interested in the behaviour of this empirical process near some fixed function $g_0$, and we denote a neighbourhood of $g_0\in\mathcal{G}$ by
\beqn
\mathcal{G}(\delta)=\{g\in\mathcal{G}:\no g-g_0\no_{2,P_X}\leq \delta\}.
\label{eq:neighbour}
\eeqn
It is important to note that the localization is performed on the set $\mathcal{G}$ since we want to estimate $g^*\in\mathcal{G}$. The aim is to measure the influence of the parameter $\lambda\in\Lambda$ in the behaviour of (\ref{eq:proc}) in the neighbourhood $\mathcal{G}(\delta)$. If we deal with a kernel deconvolution classifier, $\lambda\in\R^d$ is a set of bandwidths of a deconvolution kernel. However, in such a generality, $\lambda$ could be any kind of regularization parameter (see \cite{inversestatlearning}).

In order to apply concentration inequalities of Bernstein's type, we need the two following assumptions:
\begin{description}
\item[(A1)] for any $\lambda\in\Lambda$, there exists $b(\lambda)$ : $\sup_{g\in\mathcal{G}}\no g_\lambda \no_\infty\leq b(\lambda)$.
\item[(A2)] There exists a pseudo-distance $d$ on $\mathcal{G}$ such that $\forall g,g'\in\mathcal{G}$, $\no g_\lambda-g'_\lambda\no_{2,Z}\leq c(\lambda)d(g,g')$.
\end{description} 
These assumptions are satisfied for the particular case of the paper with Lemma \ref{lip} above where $b(\lambda)=\Pi_{i=1}^d\lambda_i^{-\beta_i-1/2}$, $d=d_\Delta$ and $c(\lambda)=\Pi_{i=1}^d\lambda_i^{-\beta_i}$. The first assumption \textbf{(A1)} is necessary to use standard uniform concentration inequalities in the bounded case (such as Bernstein or Talagrand inequalities). Moreover, \textbf{(A2)} ensures a control of the entropy of $\mathcal{G}_\lambda=\{g_\lambda,g\in\mathcal{G}\}$ thanks to standard entropy condition over the pseudo-metric space $(\mathcal{G},d)$. Indeed, using for instance \cite{wvdv}, we have under the second asumption:
\beqnn
\HH\big(\GG^\lambda,\delta,L_2(P_Z)\big)\leq\HH\left(\mathcal{G},\frac{\delta}{c(\lambda)},d\right).
\eeqnn
Next lemma proposes a control of the increments of the noisy empirical process (\ref{eq:proc}) when the class $\GG^\lambda$ satisfies the two previous assumptions. 
\begin{lemma}
\label{noisyep}
Consider a class of functions $\{g_\lambda,g\in\GG\}$ satisfying \textbf{(A1)-(A2)}. Let $g_0\in\GG$ and $\mathcal{G}(\delta)$ the set introduced in (\ref{eq:neighbour}). Suppose there exists some $0<\alpha<2$ such that:
\beqn
\label{entropyrate}
\HH_B(\GG,\delta,d)\leq c'\delta^{-\alpha}.
\eeqn
Let us consider $n_0=\inf\{n\in\mathbb{N}^*:\delta_n(\lambda)<1\}$ where
$$\delta_n(\lambda):=c(\lambda)^{\frac{\alpha}{2+\alpha}}b(\lambda)^{\frac{2}{2+\alpha}}n^{-\frac{1}{2+\alpha}}.$$
\begin{itemize}
\item Then there exist constants $c_1,c_1'$ which depend on $\alpha,c'$ such that $\forall n\geq n_0$, $\forall \delta\in[\delta_n(\lambda),1[$:
\beqnn
\P\left(\sup_{g\in\GG(\delta)}|\nu^\lambda_n(g)-\nu^\lambda_n(g^0)|\geq c_1c(\lambda)\delta^{1-\frac{\alpha}{2}}\right)\leq \exp\left(-c'_1\delta^{-\alpha}\right).
\eeqnn 
\item There exists constants $c_2,c_2'>0$ which depends on $\alpha,c'$ such that for $T\geq c_2'$, for $ n\geq n_0$:
\beqnn
\P\left(\sup_{g\notin\GG(\delta_n(\lambda))}\frac{|\nu^\lambda_n(g)-\nu^\lambda_n(g^0)|}{c(\lambda)\no g-g_0\no^{1-\frac{\alpha}{2}}}\geq T\right)\leq \exp\left(-\frac{T}{c_2}\right).
\eeqnn
\end{itemize}
\end{lemma}
The proof is an application of \cite{vdg} and consists in a noisy version of Lemma 5.13 in \cite{vdg}. This result is of pratical interest to control the variance of our deconvolution estimator. In particular, we use in the proofs of Theorem  \ref{prop:aminima} the fact that:
\beqnn
V_n=\sup_{g\in\GG}\frac{|\nu^\lambda_n(g)-\nu^\lambda_n(g^0)|}{\no g-g_0\no^{1-\alpha/2}c(\lambda)\vee c(\lambda)^{\frac{2\alpha}{2+\alpha}}n^{-\frac{2-\alpha}{2(2+\alpha)}}b(\lambda)^{\frac{2-\alpha}{2+\alpha}}}=0_\P(1),
\eeqnn
as $n\rightarrow +\infty$.\\
\textsc{Proof}
From Lemma 5.7, gathering with Lemma 5.8 of \cite{vdg}, for any $g\in\mathcal{G}(\delta)$, we have:
$$
\P(\nu_n^\lambda(g)-\nu_n^\lambda(g_0)\geq a)\leq\exp\left(-\frac{a^2}{8c(\lambda)^2\delta^2}\right),\,\forall a\leq \sqrt{n}\frac{c(\lambda)^2}{b(\lambda)}\delta^2.
$$
Next step is to use the following noisy version of Theorem 5.11. For any $a>0$ satisfying:
\beqn
\label{conditiona}
C_0\left[\sqrt{2}\delta c(\lambda)\vee \int_{a/2^6\sqrt{n}}^{\sqrt{2}\delta c(\lambda)}\sqrt{\mathcal{H}_B(\mathcal{G}^\lambda,u,L_2(P_Z))}du\right]\leq a \leq \sqrt{n}\frac{c(\lambda)^2}{K(\lambda)}\delta^2\wedge 8\sqrt{2n}\delta c(\lambda),
\eeqn
for some universal constant $C_0>0$, we have:
$$
\P(\sup_{g\in\mathcal{G}(\delta)}|\nu_n^\lambda(g)-\nu_n^\lambda(g_0)|\geq a)\leq\exp\left(-\frac{a^2}{4C\delta^2c(\lambda)^2}\right),
$$
where $C>0$ depends on $C_0$.\\
Hence from assumption \eqref{entropyrate}, for $n\geq n_0$, we have for any $\delta_n(\lambda)\leq \delta<1$, by choosing  $a=c_1c(\lambda)\delta^{1-\frac{\alpha}{2}}$ in \eqref{conditiona}:
\beqnn
\P\left(\sup_{g\in\GG(\delta)}|\nu_n^\lambda(g)-\nu_n^\lambda(g_0)|\geq c_1c(\lambda)\delta^{1-\frac{\alpha}{2}}\right)\leq \exp\left(-c'_1\delta^{-\alpha}\right).
\eeqnn 
To show the second statement, we apply the peeling device as in \cite{vdg}. Introduce:
$$
S=\inf\{s\geq 1:2^{-s}<\delta_n(\lambda)\}.
$$
Then we have, for $T=2^{1-\frac{\alpha}{2}}c_1$:
\begin{eqnarray*}
\P\left(\sup_{g\notin\mathcal{G}(\delta_n(\lambda))}\frac{|\nu_n^\lambda(g)-\nu_n^\lambda(g_0)|}{c(\lambda)\no g-g_0\no^{1-\frac{\alpha}{2}}}\geq T\right)&\leq &\sum_{s=1}^S\P\left(\sup_{2^{-s}\leq\no g-g^0\no\leq 2^{-s+1}}\frac{|\nu_n^\lambda(g)-\nu_n^\lambda(g_0)|}{\no g-g_0\no^{1-\frac{\alpha}{2}}}\geq c(\lambda)T\right)\\
& \leq & \sum_{s=1}^S\P\left(\sup_{g\in\mathcal{G}(2^{-s+1})}|\nu_n^\lambda(g)-\nu_n^\lambda(g_0)|\geq c(\lambda)c_1\left(2^{-s+1}\right)^{1-\frac{\alpha}{2}}\right)\\
&\leq &\sum_{s=1}^S\exp\left(-c'_1\left(2^{-s+1}\right)^{-\alpha}\right)=\exp\left(-\frac{T}{c_2}\right),
\end{eqnarray*}
where $c_2>0$ is a function of $\alpha,c_1$ and $c'_1$.

\begin{flushright}
$\Box$
\end{flushright}

\bibliographystyle{plain}
\bibliography{referencesbf}

\end{document}